\documentclass{article}

\usepackage{amsmath}
\usepackage{mathrsfs}
\usepackage{theorem}
\usepackage{amsfonts}

{\theorembodyfont{\slshape}
\newtheorem{theorem}{Theorem}
\newtheorem{proposition}{Proposition}
\newtheorem{lemma}{Lemma}
\newtheorem{corollary}{Corollary}

}
{\theorembodyfont{\rmfamily}

\newtheorem{remark}{Remark}

}
\def\proof{{\noindent\sc Proof. \quad}}
\def\eproof{{\mbox{}\hfill\qed}\medskip}
\newcommand\qed{{\unskip\nobreak\hfil\penalty50\hskip2em\vadjust{}
\nobreak\hfil$\Box$\parfillskip=0pt\finalhyphendemerits=0\par}}

\makeatletter
\renewcommand{\bar}{\overline}
\def\proof{{\noindent\sc Proof. \quad}}
\def\eproof{\hfill{\mbox{}\qed}}
\def\qed{\unskip\penalty10000\hbox{\enspace\vbox{\offinterlineskip
      \hrule height0.6pt\hbox{\vrule width0.6pt\vbox to 6pt{\vfil
      \hrule width 6pt height\z@}\kern-0.6pt\vrule width0.6pt}
      \kern-0.6pt\hrule height0.6pt}}\par\medbreak}
\makeatother

\newcommand{\x}{\times}

\renewcommand{\hat}{\widehat}
\renewcommand{\bar}{\overline}

\makeatletter
\def\R{\mathchoice{{\setbox0=\hbox{\rm I}\copy0\kern-0.55\wd0\hbox{\rm R}}}{%
  {\setbox0=\hbox{\rm I}\copy0\kern-0.55\wd0\hbox{\rm R}}}{%
  {\setbox0=\hbox{$\m@th\scriptstyle\rm I$}\copy0\kern-0.5\wd0%
  \scriptstyle\rm R}}{%
  {\setbox0=\hbox{$\m@th\scriptscriptstyle\rm I$}\copy0\kern-0.6\wd0%
  \scriptscriptstyle\rm R}}}
\def\N{\mathchoice{{\setbox0=\hbox{\rm I}\copy0\kern-0.55\wd0\rm N}}{%
  {\setbox0=\hbox{\rm I}\copy0\kern-0.55\wd0\rm N}}{%
  {\setbox0=\hbox{$\m@th\scriptstyle\rm I$}\copy0\kern-0.5\wd0%
  \scriptstyle\rm N}}{%
  {\setbox0=\hbox{$\m@th\scriptscriptstyle\rm I$}\copy0\kern-0.6\wd0%
  \scriptscriptstyle\rm N}}}
\def\Z{\mathchoice{{\setbox0=\hbox{\rm Z}\copy0\kern-0.5\wd0\box0}}{%
  {\setbox0=\hbox{\rm Z}\copy0\kern-0.5\wd0\box0}}{%
  {\setbox0=\hbox{$\m@th\scriptstyle\rm Z$}\copy0\kern-0.55\wd0\box0}}{%
  {\setbox0=\hbox{$\m@th\scriptscriptstyle\rm Z$}\copy0\kern-0.6\wd0\box0}}} 
\def\Q{\mathchoice{{\setbox0=\hbox{\rm Q}\kern0.2\wd0 0\kern-0.9\wd0\box0}}{%
  {\setbox0=\hbox{\rm Q}\kern0.2\wd0 0\kern-0.9\wd0\box0}}{%
  {\setbox0=\hbox{$\m@th\scriptstyle\rm Q$}\kern0.2\wd0 0\kern-0.9\wd0\box0}}{%
  {\setbox0=\hbox{$\m@th\scriptscriptstyle\rm Q$}\kern0.2\wd0 0%
  \kern-0.9\wd0\box0}}}
\def\C{\mathchoice{{\setbox0=\hbox{\rm C}\kern0.47\wd0%
  {\vrule height0.95\ht0depth1.1\dp0width0.7pt}\kern-0.47\wd0\box0}}{%
  {\setbox0=\hbox{\rm C}\kern0.47\wd0%
  {\vrule height0.95\ht0depth1.1\dp0width0.7pt}\kern-0.47\wd0\box0}}{%
  {\setbox0=\hbox{$\m@th\scriptstyle\rm C$}\kern0.47\wd0%
  {\vrule height0.95\ht0depth1.1\dp0width0.5pt}\kern-0.47\wd0\box0}}{%
  {\setbox0=\hbox{$\m@th\scriptscriptstyle\rm C$}\kern0.47\wd0%
  {\vrule height0.95\ht0depth1.1\dp0width0.4pt}\kern-0.47\wd0\box0}}} 

\def\e{{\rm e}}

\def\Cap {{\tt cap}}

\def\Oh {{\cal O}}
\def\CN{{\mathscr C}}
\def\UN{{\mathscr U}}
\def\DN{{\mathscr D}}
\def\Span {{\rm Span}\,}

\def\I {{\rm Int}\,}
\def\bfP{{\bf P}}
\def\bfE{{\bf E}}
\def\VAR{{\bf VAR}}
\def\G{{\rm G}}

\def\aff{{\rm aff}}
\def\sphere{{\rm S}}
\def\conv{{\rm conv}}
\def\AS{\mathfrak{as}}
\def\cone{{\rm cone}}
\def\T{{\rm T}}

\def\JACM{Journal of the ACM}

\def\SIOPT{SIAM Journal on Optimization}

\def\JAMS{Journal of the Amer. Math. Soc.}

\def\JoC{J. of Complexity}
\def\MP{Math. Program.}

\begin{document}

\begin{title} 
{{\bf On Tail Decay and Moment Estimates 
of a Condition Number for Random Linear Conic Systems}}
\end{title}
\author{Dennis Cheung,\quad Felipe Cucker\thanks{This author  
has been substantially funded by
a grant from the Research Grants Council of the
Hong Kong SAR (project number CityU 1085/02P).}\\
Department of Mathematics\\
City University of Hong Kong\\
83 Tat Chee Avenue, Kowloon\\
HONG KONG\\
e-mail: 
{\tt \{50003110@student,macucker@math\}.cityu.edu.hk}
\and 
Raphael Hauser\thanks{This author has been supported 
by Felipe Cucker's grant from the Research Grants Council of the
Hong Kong SAR (project number CityU 1085/02P) and 
through a grant of the Nuffield Foundation under the 
``Newly Appointed Lecturers'' grant scheme, (project number NAL/00720/G).}\\
Oxford University Computing Laboratory\\ 
Wolfson Building, Parks Road, Oxford, OX1 3QD\\ 
UNITED KINGDOM\\
e-mail: {\tt hauser@comlab.ox.ac.uk}
}

\date{}
\maketitle  

\begin{quote}
{\bf Abstract.\ } In this paper we study the distribution tails 
and the moments of $\CN(A)$ and $\log\CN(A)$, where $\CN(A)$ is 
a condition number for the linear conic system $Ax\leq 0$, $x\neq 0$, 
with $A\in\R^{n\times m}$. We consider the case where $A$ is a 
Gaussian random matrix. For this input model we characterise the 
exact decay rates of the distribution tails, we improve the existing 
moment estimates, and we prove various limit theorems for the cases 
where either $n$ or $m$ and $n$ tend to infinity. Our results are 
of complexity theoretic interest, because interior-point methods 
and relaxation methods for the solution of $Ax\leq 0$, $x\neq 0$ 
have running times that are bounded in terms of $\log\CN(A)$ and 
$\CN(A)^2$ respectively.
\end{quote}

\begin{quote}
\end{quote}

\begin{quote}
\end{quote}

\section{Introduction} 

Let $A\in\R^{n\times m}$ be given and consider the 
two systems\footnote{Usually, in the literature, 
one considers a $m\times n$ matrix $A$ appearing 
in (\ref{primalcons}), the 
``primal system'', and its transpose 
$A^{\mathrm{T}}$ appears in (\ref{dualcons}), the 
``dual system.'' We revert this notation here 
since in most of this paper we will deal with 
system (\ref{dualcons}) and we do not want to burden 
the notation with the transpose superscript.} 
\begin{equation}
Ax\leq 0,\;x\neq 0  \label{dualcons}
\end{equation}
and 
\begin{equation}
A^{\mathrm{T}}y=0,\;y\geq 0,\;y\neq0.  \label{primalcons}
\end{equation}
It is well-known that, if $A$ is full row rank, 
one of these systems has a strict solution (one
for which the satisfied inequality is strict in all coordinates) 
if and only if the other has no solutions at all.

The following feasibility problem is a standard subproblem in 
linear programming:
\begin{quote}
{\bf (FP)}  Given a $n\times m$ full row rank real matrix $A$, 
  decide which of
  (\ref{dualcons}) or (\ref{primalcons}) is strictly 
  feasible and return a strict solution for it.
\end{quote}
Iterative algorithms that solve this problem, such as variants of 
interior-point or ellipsoid methods, have a cost which depends on some 
measure of conditioning of the matrix $A$ besides the natural 
dependence on $n$ and $m$. For instance, 
a finite-precision algorithm solving the problem 
above is analysed 
in~\cite{CP01} to show a complexity of  
\begin{equation}\label{javier}
    \Oh\left((m+n)^{3.5}\left(\log(m+n)+
     \log C(A)\right)^3\right).
\end{equation}
Here $C(A)$ may be either the condition number introduced by Renegar 
in~\cite{Renegar94b,Renegar95,Renegar95b}, which we denote by 
$C_R(A)$, or a generalisation of Goffin's ``inner measure'' 
\cite{Goffin1,Goffin2} introduced in~\cite{ChC00}, which we denote by 
$\CN(A)$ in the sequel.

Another family of algorithms whose complexity can be analysed in 
terms of $\CN(A)$ are the many variants of the
Agmon-Motzkin-Sch\"onberg (AMS) relaxation method
\cite{Agmon,Motzkin} for solving systems of linear inequalities.
This includes the cyclic projection method, the perceptron algorithm 
and certain types of subgradient algorithms. The complexity of these 
methods is typically proportional to $\CN(A)^2$. 
For example, for feasible systems the perceptron algorithm is 
guaranteed to find a solution to $Ax<0$ in 
\begin{equation}\label{perceptron}
\Oh\left(\CN(A)^2\right)
\end{equation}
iterations (see Appendix~B for further remarks). 
Although less interesting from a complexity perspective, such 
algorithms have appealing aspects in applications where $m$ 
and $n$ are both very large, for example in 
tumour radiation therapy planning. AMS relaxation is also the historic 
context in which the condition number $\CN(A)$ has first been studied, 
albeit only in the case of feasible systems, see
\cite{Goffin1,Goffin2}. 

It is also worthwhile mentioning that there exist close links between 
the condition number $\CN(A)$ and the notion of {\em margin} that
plays a key role in the learning theory literature. Furthermore, 
the condition number $\CN(A)$ has applications in the backward error 
analysis and in estimating the stability of linear feasibility problems. 
For all these reasons, studying the moments of both $\log\CN(A)$ and 
$\CN(A)$ will be interesting in the probabilistic setting outlined 
below. 

Unlike $m$ and $n$, the condition number of $A$ is not immediate 
to determine from the data $A$ and seems to require a 
computation which is not easier than solving the feasibility problem 
instance described by $A$ (see~\cite{Renegar94} for a discussion). 
In addition, there are no bounds on its magnitude as a function of $m$ 
and $n$. It may actually be infinite. Thus, bounds such as 
\eqref{javier} and \eqref{perceptron} tell us little about the running 
time we can expect for a given input $A$.

A reasonable way to cope with this situation is to assume a
probability measure on the space of $n\x m$ 
matrices $A$ and to estimate the expected value of 
$\log\CN(A)$ (or that of $\CN(A)$). A standard choice of distribution is the 
Gaussian model. We say that a random matrix is {\em Gaussian} 
when its entries are i.i.d.~random variables with standard 
normal distribution. The main result in~\cite{ChC01} shows that if 
$A$ is a Gaussian $n\times m$ matrix,  then 
\begin{equation}\label{thChC}
\bfE\left[\log\CN(A)\right]=
\left\{
\begin{array}{ll}
\mathcal{O}(\min \{m\log n,n\}) & \mbox{if $n>m$} \\[8pt] 
\mathcal{O}(\log n) & \mbox{otherwise}
\end{array}\right.
\end{equation}
and
$$
  \bfE[\log C_R(A)]\leq\bfE[\log\CN(A)]+\frac{5\log n}{2}
   +\frac{\log m}{2}+2\log 2.
$$
\hspace{1cm}

In most practical occurrences of the feasibility 
problem {\bf (FP)} 
one deals with matrices for which $n$ is some orders of 
magnitude larger than $m$. The case $n\gg m$ ($n$ much 
larger than $m$) is actually the case of interest 
among researchers in linear programming. 
In~\cite{CW01} the estimate \eqref{thChC} was refined 
to prove that, when $n$ is moderately larger than $m$ one has   
$\bfE[\log\CN(A)]\leq \max\{\log m,\log\log n\}+\Oh(1)$. 

The main contributions of this paper are the following:

\noindent
{\bf (i)}\quad
We further strengthen the existing bounds on $\bfE[\log\CN(A)]$: 
In Corollary \ref{polycor2} we show that if $n\geq m$
\begin{equation*}
\bfE[\log\CN(A)]\leq m\log 2
\end{equation*}
asymptotically as $m\rightarrow\infty$, and in Corollary \ref{polycor3} 
we show that if $n\geq 5m$ and $\gamma>0$ is an arbitrary constant, 
then 
\begin{equation*}
\bfE[\log\CN(A)]\leq m^{\gamma}
\end{equation*}
asymptotically as $m\rightarrow\infty$. 

\noindent
{\bf (ii)}\quad
We generalise the bounds to arbitrary moments of $\log\CN(A)$ and 
derive similar bounds for the moments of $\CN(A)$: see Corollaries 
\ref{estimates},\ref{polycor} and \ref{CORI}. 

\noindent
{\bf (iii)}\quad
We derive various limit theorems that are of interest in large-scale 
problems: in the case where $m$ is fixed and $n\rightarrow\infty$, 
we show in particular that 
\begin{align*}
\CN(A)&\stackrel{n\rightarrow\infty}{\longrightarrow}_P 1,\\
\lim_{n\rightarrow\infty}\bfE[(\log\CN(A))^{\gamma}]&=0\quad
\forall\gamma>0,\\
\lim_{n\rightarrow\infty}\bfE[(\CN(A))^{\gamma}]&=1\quad\forall 
\gamma\in[0,1),
\end{align*}
see Corollaries \ref{corollary 4},\ref{limits} and \ref{corollary 6}. 
In the case where $n\geq m$ and 
$m\rightarrow\infty$, it is again Corollaries \ref{polycor2} and 
\ref{polycor3} that bound the asymptotic growth rates of 
$\bfE[\log\CN(A)]$. 

\noindent
{\bf (iv)}\quad  
Previous probabilistic analyses of $\log\CN(A)$ relied on the fact 
that linear algebraic operations applied to Gaussian matrices
(including vectors as a special case) lead again to Gaussian
matrices. The analysis presented here is very different because it 
is based on geometry on high dimensional spheres. This approach 
makes it possible to derive not just upper bounds on the tail 
decay of $\log\CN(A)$, but also lower bounds: Theorems \ref{tail_th} 
and \ref{lowbd} show that there exist functions $c(m,n)\geq d(m,n)$ 
of $m$ and $n$ such that 
\begin{equation*}
\frac{d(m,n)}{t}\leq\bfP\left[\CN(A)\geq t\right]\leq\frac{c(m,n)}{t}
\end{equation*}
for all $t$ large enough. This implies that the distribution tails of 
$\log\CN(A)$ asymptotically decay exactly at the exponential 
rate $\e^{-t}$, see Corollary \ref{exact decay rates}. 
More importantly, the geometric 
analysis we developed here generalises to almost arbitrary probability 
measures that are absolutely continuous with respect to the uniform 
measure on the sphere, as we will show in a follow-up paper. In 
the general case, the tail decay rates of $\log\CN(A)$ are again 
exponential, but the exponent depends on a parameter defined as 
a function of the distribution. This exponential decay is perhaps 
the most important conclusion of our analysis, as it explains why 
the polyhedral feasibility problem -- and linear programming by 
extension -- is ``empirically strongly polynomial''. We discuss the 
relevance of this notion in Section \ref{discussion}. 

\section{Complexity Theoretic Context}\label{discussion} 

\subsection{Background}

The interest in the distribution of $\log\CN(A)$ stems to a large 
extent from the conjecture that there exist so-called {\em strongly 
polynomial} algorithms for linear programming. Let us give a brief 
explanation for the uninitiated reader: since the mid-1940s, 
variants of Dantzig's simplex method proved to be efficient algorithms 
for solving linear programming problems in practice, despite the 
fact that in the worst case these methods terminate only after a 
number of iterations that is exponential in the ``size'', or the total 
input data length of the problem. As interest in complexity theory
grew, many researchers believed that a good algorithm should terminate 
within a number of iterations that is bounded by a polynomial in the input 
size. Thus, the simplex method is not a polynomial algorithm.
 
Surprising new approaches to linear programming subsequently proved to 
be polynomial time algorithms for linear programming under the Turing 
model: Khachiyan's ellipsoid method~\cite{Kat79}, Karmarkar's
method \cite{Karmar84} and the many interior-point methods
developed since then are algorithms of this kind. These algorithms 
are guaranteed to terminate in polynomial time when the input data 
are rational and the input size is measured by the total bit-length 
of the data. 

On the other hand, under the so-called {\em real} complexity 
model one considers linear programming problems whose input data 
are real numbers and imagines a hypothetical computer that can perform 
operations on real numbers. In this model the complexity of an
algorithm is the number of operations that are needed in the worst 
case to solve the problem. Such an algorithm is called {\em strongly 
polynomial} if its complexity is a polynomial function in the number 
of constraints and variables (the ``dimension'') of the underlying 
problem. Neither the ellipsoid method nor any of the known interior-point 
algorithms for linear programming is known to be strongly polynomial. 
In fact, their running time is theoretically unbounded! This 
is in stark contrast to the simplex method which is guaranteed to 
terminate in exponential time. 

The situation is not hopeless, however, for the real complexity 
of ellipsoid and interior-point methods can be bounded by a 
polynomial in the problem dimension and the logarithm of a condition 
number, see the results cited in the introduction. Earlier 
relevant papers on this subject and on other applications of 
LP condition numbers include (among others) 
\cite{Renegar94b,Renegar95,
Renegar95b,FreundVera,FreundVera99,Vera98,CP01,VavasisYe,VavYe2}. 

This condition-based complexity analysis is not new. 
In numerical linear algebra it occurs, for instance, 
in the analysis of the conjugate gradient method 
(cf.~\cite{Reid,TrefethenBau}); in polynomial 
equation solving it occurs in the analysis of homotopy 
methods~\cite{Bez1} or of grid-based 
methods~\cite{CS98,Cucker99b}. A recent survey for 
linear programming is~\cite{ChCYe}. A conceptually related idea in 
discrete mathematics is that of 
{\em parameterized complexity}~\cite{DowFe}. 

The question of whether linear 
programming is strongly polynomial time solvable is considered an 
important open problem and has many ramifications within complexity 
theory; in his list of 18 mathematical problems for the XXI 
century~\cite{Smale98} Steve Smale includes this question as 
Problem~9. 

An interesting approach to get around the difficult issue of strong 
polynomiality of linear programming is the average case analysis of 
algorithms. The average case analysis reveals that linear programming 
is strongly polynomial time solvable {\em on average}, that is, 
there exist algorithms whose average running times are polynomial 
under the real model when the input data are normally distributed. 
The simplex method was known to possess this property since the 
early 1980s~\cite{Borgwardt1,Borgwardt2,Borgwardt3,Smale83}. 
Similar work was continued in \cite{AdlerMegiddo,Megiddo86,Todd86} 
and \cite{Borgwardt4}. More recently, 
the attention has shifted to the average case analysis of
interior-point methods \cite{Anstreicher,Todd,Minoux,Huhn-Borgwardt-I} 
and \cite{Huhn-Borgwardt-II}. While all of these papers followed an 
analysis 
pertaining to particular algorithms, it is also possible to derive 
similar results by analysing the expected value of condition numbers 
under Gaussian (or other) input data. Combined with a 
condition-based complexity analysis this yields average case running time 
bounds for particular algorithms. This was the approach pursued in 
\cite{TTY98,ChC01} and \cite{CW01}, and it is also the approach we 
pursue in the present paper. 

We should point out that the relevance of average case analyses is 
subject to some justified scrutiny which we will further address and 
respond to in the next paragraph.

\subsection{Discussion: Relevance of Our Results}

As mentioned in the synopsis of the introduction, the results we 
will derive in this paper include in particular polynomial bounds 
of $\bfE[\log\CN(A)]$ in $m$ and $n$. In the literature on the 
probabilistic analysis of linear programming such 
results are considered important because they show that LP is
``strongly polynomial on average''. Precisely how significant are 
such statements from a complexity theoretic view point?

On the one hand, the average behaviour of an algorithm on 
random input data yields in itself an interesting complexity measure 
which, as many would argue, can be more relevant than the study of the 
worst case scenario. The weakness of this argument is that the relevance 
of average case results depends on how well the assumed probabilistic 
model describes the distribution of the input data that one might 
observe in particular applications. Without doubt, uniform or Gaussian 
matrices are an inadequate model in most but a few applications. In 
a follow-up paper we will therefore show how the techniques developed 
here extend to matrices with much more general probability
distributions. We chose to treat the Gaussian case separately because 
it can be presented in a non-measure-theoretic setting which is 
accessible to a wider audience, and because it allows to directly 
compare our results, which were obtained by arguments based on 
spherical geometry, with the results obtained in earlier papers 
via the very different techniques of linear algebra on Gaussian
matrices. 

On the other hand, it is sometimes argued that understanding the 
average behaviour of interior-point algorithms constitutes a step 
towards proving strong polynomiality of linear programming. This 
argument is somewhat weaker than the first, because it may of course 
be that linear programming is ``strongly polynomial on average'' 
for a wide range of input distributions whilst not allowing a 
strongly polynomial time algorithm. Thus, the two phenomena 
might be unrelated.

In our view the most relevant link between the results presented 
in this paper and the conjectured strong polynomiality of linear 
programming (and the closely associated linear feasibility problem 
treated here) consists not in 
bounds on the moments of $\log\CN(A)$ but in the exponential decay 
of its distribution tails observed in Corollary~\ref{exact decay
rates}. This fact explains why algorithmic experiments have a tendency 
to strengthen the intuition that the conjecture of strong
polynomiality be true: $\Oh(e^t)$ simulations are needed to observe 
an event in which $\log\CN(A)>t$ (and even in that case it is not 
guaranteed that an algorithm is necessarily slow at solving the
problem). Thus, it is impossible to observe the really bad cases
in random experiments. In contrast, in algorithms whose complexity 
depends polynomially on $\CN(A)$, much fewer simulations reveal 
cases with extreme running times, because it takes $\Oh(t)$ 
experiments to detect an event in which $\CN(A)>t$. Of course, this 
argument is again subject to the criticism that the exponential 
decay of $\bfP[\CN(A)>t]$ might be a particularity of the chosen 
input distribution for the data of $A$. However, our analysis of 
more general distributions shows that exponential decay rates are 
a common feature of a very general family of distributions.

These observations suggest a notion of ``empirical strong 
polynomiality'': if linear programming is not strongly polynomial 
time solvable, then this fact cannot be observed in 
random experiments because the relevant events are exponentially
rare. This is in our view the correct interpretation of our results 
and the main message of this paper.

\section{Basic Notions and Notation}\label{notions}

If there exists $x\in\R^m$, $x\neq0$, 
such that $Ax\leq 0$ we say that $A$ 
is {\em feasible}. Otherwise, we say that $A$ is 
{\em infeasible}. Also, if there exists a vector $x$ such that 
$Ax<0$ componentwise, then we say that $A$ is {\em strictly
feasible}. If $A$ is feasible but not strictly feasible then we 
say that $A$ is {\em ill-posed}.

In the sequel we consider $\arccos(t)$ as a function from $[-1,1]$
into $[0,\pi]$. In this region, both $\cos$ and $\arccos$
are decreasing functions.

Let $a_i$ be the $i$th row of $A$. Denote by 
$\theta_i(A,x)$ the angle between $a_i$ and $x$, 
that is, $\arccos(\frac{a_i\cdot x}{\|a_i\|\|x\|})$. Let 
$\theta(A,x)=\min_{1\leq i\leq n}\theta_i(A,x)$
and $\bar x$ be any vector in $\R^m\setminus\{0\}$, s.t.
$\theta(A)=\theta(A,\bar x)=\sup_{x\in\R^m}\theta(A,x)$.
The condition number $\CN(A)$ is defined as
\begin{equation}\label{cond}
\CN(A)=|\cos(\theta(A))|^{-1}.
\end{equation}

In the case where the system $Ax\leq 0$ is feasible, $\CN(A)$ 
is the same as Goffin's condition number $\nu$ \cite{Goffin1,
Goffin2}, which he developed to analyse step length rules that 
guarantee finite convergence of the relaxation method applied 
to feasible systems of linear inequalities. $\CN(A)$ is also 
closely related to Agmon's condition number $\lambda$ \cite{Agmon} 
with which it coincides in some cases, but again $\lambda$ is 
only defined for feasible systems. Thus, $\CN(A)$ is more 
general because it is defined both for feasible and infeasible 
systems. 

It is not difficult to see that $A$ is strictly feasible iff 
$\theta(A)>\pi/2$, ill-posed iff $\theta(A)=\pi/2$ and infeasible 
iff $\theta(A)<\pi/2$. It is also easy to show, 
using a compactness argument, that a vector $\bar{x}$ such 
as used in the definition of $\CN(A)$ exists. Note that since 
$\CN(A)$ is defined purely in terms of angles between vectors, 
$\CN$ is invariant under positive rescaling of the rows of $A$. 
Hence, we may assume without loss of generality that all 
rows of $A$ have been rescaled to unit length. Our analysis then 
reduces to geometry on the unit sphere.

Let $\sphere^{m-1}$ denote the unit sphere in $\R^m$. 
For $p\in \sphere^{m-1}$ and $\rho\in[0,\pi]$ we denote 
by $\Cap(p,\rho)$ the circular cap with centre in $p$ and angular 
radius $\rho$, that is,  
$$
 \Cap(p,\rho)=\left\{x\in \sphere^{m-1}:x\cdot p\geq \cos\rho\right\}. 
$$
By $\partial\Cap(p,\rho)$ and $\I(\Cap(p,\rho))$ we denote 
the boundary and interior of $\Cap(p,\rho)$, respectively in the 
standard topology on $\sphere$, that is, 
\begin{align*}
\partial\Cap(p,\rho)&=\left\{x\in \sphere^{m-1}:x\cdot p=
\cos\rho\right\},\\ 
\I(\Cap(p,\rho))&=\{x\in \sphere^{m-1}:x\cdot p>\cos\rho\}. 
\end{align*}
\hspace{1cm}

The following result provides geometric insight about the condition 
number $\CN(A)$. 

\begin{proposition}\label{geo1}
It is true that $0<\theta(A)\leq\rho\leq\pi$ if and only if 
$\bigcup_{i=1}^n\Cap(a_i,\rho)=\sphere^{m-1}$.
\end{proposition}

\proof
\begin{eqnarray*}
\theta(A)\leq\rho
&\Leftrightarrow&\max_{x\in\R^{m}\setminus\{0\}}\theta(A,x)\leq\rho\\
&\Leftrightarrow& \forall x\in \sphere^{m-1}, 
  \min_{1\leq i\leq n}\theta_i(A,x)\leq\rho\\
&\Leftrightarrow& \forall x\in \sphere^{m-1}, 
  \exists i\in\{1,2,\ldots,n\}, x\in\Cap(a_i,\rho)\\
&\Leftrightarrow&\bigcup_{i=1}^n\Cap(a_i,\rho)=\sphere^{m-1}.
\end{eqnarray*}
\eproof
\hspace{1cm}

\section{Characterising Extremal Circular Caps}\label{caps}

In this section $A$ is a real $n\times m$ matrix with unit row vectors 
$a_i$. A {\em largest circular cap} excluding rows of $A$ (in short, a  
LCP of $A$) is a cap 
$\Cap(p^*,\rho^*)$ corresponding to a maximiser $p^*$ and 
its objective function value $\rho^*=\rho(p^*)$ for the optimisation 
problem 
\begin{equation}\label{maxi}
\max_{p\in\sphere^{m-1}}\rho(p),
\end{equation}
where $\rho(p)=\min_{i}\arccos(p\cdot a_i)$. A {\em smallest circular 
cap} containing all of the $a_i$ $(i=1,\dots,n)$ (in short, a SCP of 
$A$) is a complement of a LCP of $A$.

The following result explains why the notions of LCP and SCP are 
important in the analysis of the condition number $\CN(A)$.

\begin{proposition}\label{geo2}\label{geo3}
Let $\Cap(p,\rho)$ be a LCP of $A$. Then $x^*=p$ maximises the function 
$\theta(A,x)$ and it is true that $\theta(A)=\theta(A,p)=\rho$, that 
is, $\CN(A)=1/|\cos\rho|$.
\end{proposition}

\proof 
Let $\bar{x}\in\sphere^{m-1}$ be a maximiser of $\theta(A,x)$. Then 
$\theta(A,\bar{x})=\theta(A)$ and $\I(\Cap(\bar x,\theta(A)))$ 
contains none of the rows of $A$. This shows that 
$\theta(A)\leq\rho$. On the other hand, since 
$a_i\in\Cap(-p,\pi-\rho)$ for all $i$, we have 
\begin{equation*}
\rho\leq\min_{1\leq i\leq n}\arccos(p\cdot a_i)=
\theta(A,p)\leq\theta(A).
\end{equation*}
\eproof
\hspace{1cm}

It is worth investigating the properties of LCPs and SCPs a bit 
further. A simple compactness argument shows that LCPs and SCPs exist 
for any $A$. If $A$ is 
infeasible then these caps are not unique in general. To visualise 
this fact, it is helpful to consider the limiting case where a 
countable set of points densely fills what is left of the unit 
sphere after two circular caps of equal radius but different 
midpoints have been removed.

If $A$ is strictly feasible however, 
there exists a unique LCP, respectively SCP, as the following 
convexity argument shows: it is easy to see that the strict 
feasibility of $A$ implies that $p^*\cdot a_i>0$ for all $i$ and 
$0\leq\varrho^*<\pi$ for any SCP $\Cap(p^*,\varrho^*)$. It 
suffices therefore to argue that the minimisation problem 
\begin{equation}\label{mini}
\begin{split}
\min_{p\in\sphere^{m-1}}&\varrho(p)\\
\text{s.t. }&p\cdot a_i>0\quad(i=1,\dots,n)
\end{split}
\end{equation}
has a unique local minimiser, where $\varrho(p)=\pi-\rho(-p)=
\max_{i}\arccos(p\cdot a_i)$. Suppose $p_1\neq p_2$ are two distinct 
local minimisers of \eqref{mini} such that 
$\varrho_2:=\varrho(p_2)\geq\varrho_1:=\varrho(p_1)$. Then 
\begin{equation}\label{varrho}
p_j\cdot a_i\geq\cos\varrho_j>0\quad(i=1,\dots,n; j=1,2).
\end{equation}
For $\lambda\in(0,1)$ let $p(\lambda):=\lambda p_1+
(1-\lambda)p_2$ and $\hat{p}(\lambda):=p(\lambda)/\|p(\lambda)\|$. 
Then $p_1\neq p_2$ implies that $\|p(\lambda)\|<1$ and 
\begin{equation}\label{One}
\hat{p}(\lambda)\cdot a_i>p(\lambda)\cdot a_i
\stackrel{\text{\eqref{varrho}}}{\geq}\lambda\cos\varrho_1+
(1-\lambda)\cos\varrho_2\quad(i=1,\dots,n).
\end{equation}
It follows from \eqref{One} that $p_2$ is not a local minimiser 
of \eqref{mini} if $\varrho_2>\varrho_1$, that is, it must be 
true that $\varrho_1=\varrho_2$. But then, for all $\lambda\in(0,1)$, 
\eqref{One} shows that $\Cap\bigl(\hat{p}(\lambda),\varrho(\lambda)
\bigr)$ contains all $a_i$ $(i=1,\dots,n)$, where $\varrho(\lambda)=
\arccos\bigl(\min_{i}\hat{p}(\lambda)\cdot a_i\bigr)<\varrho_1$, 
contradicting the local optimality of 
$p_1$ and $p_2$. This shows that the SCP and, by extension, the 
LCP of $A$ are unique, as claimed.

Next we investigate the idea of blocking sets. Let $\Cap(p,\rho)$ 
be a LCP of $A$. We say that 
\begin{equation}\label{proof of support1}
S=\bigl\{i:a_i\in\partial\Cap(p,\rho)\bigr\}
\end{equation}
is the {\em blocking set} of $\Cap(p,\rho)$. The blocking set corresponds to
the vectors that locally keep the LCP from having a larger 
radius.\footnote{We say ``locally'' in the sense that the largest cap 
not containing any of the $a_i$ and centred at a point $q$ has radius 
smaller than $\rho$ for all $q$ in a small neighbourhood around 
$p$. Note, however, that the blocking set does not prevent a
point $q$ far away from $p$ from being the centre of an even
larger cap not containing any of the $a_i$. This idea of a
blocking set thus explains why $\Cap(p,\rho)$ is a {\em local}
minimiser of problem (\ref{mini}). Of course, a LCP is defined as
a {\em global} minimiser of this problem. Thus, the
existence of a blocking set is a local optimality condition. 
As always in nonlinear optimisation, useful global optimisation
criteria don't really exist.} 
In fact, $S$ is the active set in the following equivalent
reformulation of \eqref{maxi}:
\begin{align*}
\max_{p\in\sphere^{m-1}}&\rho\\
\text{s.t. }&p\cdot a_i\leq\cos\rho\quad(i=1,\dots,n).
\end{align*}
By straightforward extension we also speak of the blocking set of 
a SCP.

If $A$ is strictly feasible then the blocking set of the 
(unique) LCP can have any 
cardinality $\geq\min(2,n)$, as can easily be illustrated by the 
example of three points on $\sphere^2$ that lie on a single grand 
circle and within a sufficiently small angle of one another.

The situation is rather different if $A$ is infeasible and $n\geq m$. 
In this case all blocking sets have cardinality $\geq m$. In fact, 
let $\Cap(p,\rho)$ be a LCP of an infeasible $A$, that is, 
$\theta(A)=\rho<\pi/2$, let $S$ be the blocking set of 
$\Cap(p,\rho)$ and suppose that $|S|< m$. Then there exists a vector 
$u\in \sphere^{m-1}\cap\Span(S)^{\perp}$ such that $u\cdot p\geq
0$. Let $p_{\delta}=p+\delta u$. For $\delta>0$ we have 
$\|p_{\delta}\|^2=1+\delta^2+2\delta u\cdot p>1$. Let 
$\hat{p}_{\delta}=p_{\delta}/\|p_{\delta}\|$ and $\rho_{\delta}
=\arccos\bigl(\cos\rho/\|p_{\delta}\|\bigr)$. Then $\rho_{\delta}>
\rho$. For all $i\notin S$ and $\delta>0$ sufficiently small, 
$a_i\notin\Cap(\hat{p}_{\delta},\rho_{\delta})$, since $a_i\notin 
\Cap(p,\rho)$ and $\Cap(\hat{p}_{\delta},\rho_{\delta})$ varies 
continuously as a function of $\delta$. Moreover, for $i\in S$ 
we have $p_{\delta}\cdot a_i=p\cdot a_i+\delta u\cdot a_i=
p\cdot a_i=\cos\rho$. Therefore $\hat{p}_{\delta}\cdot a_i
=\cos\bigl(\rho/\|p_{\delta}\|\bigr)=\cos\rho_{\delta}$, and this shows that 
$a_i\in\partial\Cap(p_{\delta},\rho_{\delta})$ for all $i\in S$.
We conclude that $\I(\Cap(\hat{p}_{\delta},\rho_{\delta}))$ contains 
none of the $a_i$, and since $\rho_{\delta}>\rho=\theta(A)$ this 
is a contradiction. Therefore, our assumption was wrong and $|S|\geq
m$.

Let us summarise what we have found so far.

\begin{proposition}\label{support1}
Let $A\in\R^{n\times m}$ have unit rows. Then 
\begin{itemize}
\item[(i)\;] If $A$ is strictly feasible there exists a unique LCP 
and SCP of $A$ but the cardinality of the blocking set is arbitrary. 
\item[(ii)\;] If $A$ is infeasible then there exist LCPs and SCPs of 
$A$ which are not unique in general, but their blocking sets always 
have cardinality $\geq m$. 
\end{itemize}
\end{proposition}

\proof
See arguments above.
\eproof
\hspace{1cm}

Let us further explore the link between blocking sets and extremal 
circular caps. We consider the set of index sets of cardinality $m$, 
${\mathcal P}_m=\{S\subset\{1,\dots,n\}: |S|=m\}$. 
If $S\in{\mathcal P}_m$ and $A\in\R^{n\times m}$ we denote by $A_S$ 
the $m\times m$ matrix obtained by removing 
all rows from $A$ with index not in $S$. Let $e=(1\dots
1)^{\T}\in\R^m$. If $A$ is nonsingular, the vectors 
\begin{equation}\label{u_S}
u_S=A_S^{-1}e\quad\text{ and }\quad
\hat{u}_S=\frac{u_S}{\|u_S\|}
\end{equation}
are well defined. 

\begin{lemma}\label{support5}
Let $S\in{\mathcal P}_m$ be such that $A_S$ is nonsingular, and 
let $p\in\sphere^{m-1}$, $\rho\in[0,\pi/2)$ and 
$S\in{\mathcal P}_m$ be such that $a_i\in\partial\Cap(p,\rho)$ 
for $(i\in S)$. Then $p=\hat{u}_S$ and $\cos\rho=\|u_S\|^{-1}$.
\end{lemma}

\proof
If $a_i\in\partial\Cap(p,\rho)$ then $a_i\cdot p=\cos\rho$.
Therefore, $A_Sp=(\cos\rho) e$ and $p=(\cos\rho)A_S^{-1}e=
(\cos\rho)u_S$. Since $\|p\|=1$, we have $|\cos\rho|=\|u_S\|^{-1}$ 
and, using $\rho<\pi/2$, $\cos\rho=\|u_S\|^{-1}$. 
We conclude that $p=\hat{u}_S$. 
\eproof
\hspace{1cm}

Blocking sets are the key tool that allow us to gain information about 
the distribution of $\CN(A)$ when the unit rows of $A$ are random
vectors with known distribution. The fact that the blocking set $S$ 
may have cardinality $|S|<m$ for strictly feasible $A$ is an obstacle 
to this analysis that requires further attention. The following two 
lemmas allow to overcome this problem in the analysis of upper and 
lower bounds on the tails of $\CN(A)$ respectively.

\begin{lemma}\label{support3}
Let $A\in\R^{n\times m}$ have $n\geq m$ unit rows, and let 
$\Cap(p,\rho)$ with $\rho<\pi/2$ contain all rows of $A$. 
Then there exist $p'\in\sphere^{m-1}$ 
and $\rho\leq\rho'<\frac{\pi}{2}$ such that $\Cap(p',\rho')$ also 
contains all rows of $A$ and at least $m$ of them lie on the 
boundary $\partial\Cap(p',\rho')$. 
\end{lemma}

\proof
Let $S=\bigl\{i:a_i\in\partial\Cap(p,r)\bigr\}$, and 
let us assume that $|S|<m$. Then there exists a vector 
$u\in\sphere^{m-1}\cap\Span(S)^{\perp}$, where we set 
$\Span(S)^{\perp}=\R^m$ if $S=\emptyset$, such that there exists 
an index $i\notin S$ with $a_i\cdot u\neq 0$. Without loss of 
generality we may assume that $p\cdot u\geq 0$. Let $i^*$ be a 
minimiser of 
$\min_{i\not\in S}\left|\frac{\cos\rho-a_i\cdot p}{a_i\cdot u}
\right|$ and 
$$
\delta=\frac{\cos\rho-a_{i^*}\cdot p}{a_{i^*}\cdot u}.
$$
Let finally $p_{\delta}=p+\delta u$. Then $p\cdot u\geq 0$ implies 
that $\|p_{\delta}\|\geq 1$ and hence, $\hat{p}_{\delta}=\|p_{\delta}\|
^{-1}p_{\delta}$ and $\rho_{\delta}=\arccos\bigl(\|p_{\delta}\|
^{-1}\cos\rho\bigr)\in[\rho,\frac{\pi}{2})$ are well defined. 
The proof of Proposition~\ref{support1} 
shows that $a_i\in\partial\Cap(\hat{p}_{\delta},\rho_{\delta})$ 
for all $i\in S$. Moreover, the definition of $\delta$ shows that 
$a_{i*}\cdot p_{\delta}=\cos\rho$ and hence that $a_{i*}\cdot 
\hat{p}_{\delta}=\cos\rho_{\delta}$. This shows that $a_{i^*}\in
\partial\Cap(\hat{p}_{\delta},\rho_{\delta})$ and thus, $|S'|>|S|$, 
where 
\begin{equation}\label{proof of support3}
S'=\bigl\{i:a_i\in\partial\Cap(\hat{p}_
{\delta},\rho_{\delta})\bigr\}. 
\end{equation}
Note also that if $a_i\notin\Cap(\hat{p}_{\delta},\rho_{\delta})$ then 
$i\notin S$ and $a_{i}\cdot\hat{p}_{\delta}<\cos\rho_{\delta}$, that is, 
$a_{i}\cdot p_{\delta}<\cos\rho$ which contradicts the choice of 
$i^*$. Therefore, $\Cap(\hat{p}_{\delta},\rho_{\delta})$ contains 
all rows of $A$. Using this construction recursively, we eventually 
arrive at the desired $p'$ and $\rho'$. 
\eproof
\hspace{1cm}

\begin{lemma}\label{support4}
Let $a_1,\dots,a_m$ be linearly independent elements of 
$\partial\Cap(p,\rho)\subset\sphere^{m-1}$ with $\rho\in[0,\pi/2)$. 
Let $\pi_p:\R^m\rightarrow\Span(p)^{\perp}$ be the orthogonal 
projection along $p$. Then $\Cap(p,\rho)$ is the SCP of 
$A=[a_1\dots a_m]^{\T}$ if and only if $0\in\conv(\pi_p a_1,
\dots,\pi_p a_m)$. 
\end{lemma}

\proof
Let us first remark that is trivial to see that $\rho<\pi/2$ 
implies that 
\begin{equation}\label{trivial}
0\in\conv(\pi_p a_1,\dots,\pi_p a_m)
\Leftrightarrow p\in\cone(a_1,\dots,a_m).
\end{equation} 
Note also that the linear independence of the $a_i$ implies that 
$A$ is strictly feasible and the SCP of $A$ is unique.
 
We will now show the ``only if'' part of the lemma. Suppose $0\notin
\conv(\pi_p a_1,\cdots,\pi_p a_m)$. Because of \eqref{trivial},  
Farkas' lemma implies that there exists $q\in\sphere^{m-1}$ such 
that $q\cdot a_i>0$ $(i=1,\dots,m)$ and $q\cdot p<0$. For 
$\delta>0$ let 
\begin{equation}\label{delta}
p_{\delta}=p+\delta q,\quad 
\hat{p}_{\delta}=\frac{p_{\delta}}{\|p_{\delta}\|}\quad
\text{ and }\quad
\rho_{\delta}=\arccos\min_{i}\bigl(a_i\cdot\hat{p}_{\delta}\bigr).
\end{equation} 
Then $\|p_{\delta}\|^{-1}=1-\delta p\cdot q+\Oh(\delta^2)$ and 
\begin{equation*}
a_i\cdot\hat{p}_{\delta}\geq
(\cos\rho+\delta a_i\cdot q)(1-\delta p\cdot q+\Oh(\delta^2))>
\cos\rho\qquad(i=1,\dots,m)
\end{equation*}
for all $\delta>0$ small enough. Therefore, $\Cap(\hat{p}_{\delta},
\rho_{\delta})$ contains all $a_i$ and $\rho_{\delta}<\rho$ for 
$0<\delta\ll 1$, showing that $\Cap(p,\rho)$ is not the SCP of $A$. 

It remains to show the ``if'' part of the lemma. Because of
\eqref{trivial} we may assume that $p=\sum_{i=1}^m\lambda_i a_i$ 
for some $\lambda_1,\dots,\lambda_m\geq 0$. Let us assume that 
$\Cap(p,\rho)$ is not the SCP of $A$. 
Then by a construction similar to \eqref{One} 
there exists a direction $d\in\Span(p)^{\perp}$ such that 
$\Cap(\hat{p}_{\delta},\rho_{\delta})$ contains all $a_i$ and 
$\rho_{\delta}<\rho$ for $0<\delta\ll 1$, where $\hat{p}_{\delta}$ 
and $\rho_{\delta}$ are defined as in \eqref{delta}. That is 
to say, 
\begin{equation*}
a_i\cdot\hat{p}_{\delta}=\frac{a_i\cdot p+\delta a_i\cdot d}
{\sqrt{1+\delta^2}}>\cos\rho=a_i\cdot p
\end{equation*}
for $(i=1,\dots,m)$ and $0<\delta\ll 1$, which shows that 
$a_i\cdot d>0$ for all $i$. But then 
$0=p\cdot d=\sum_{i=1}^m\lambda_i a_i\cdot d>0$, which shows that 
our assumption was wrong and $\Cap(p,\rho)$ is indeed the SCP 
of $A$. 
\eproof
\hspace{1cm}

\section{The Input Distribution}\label{probability model}

It is well known that if $A$ is a Gaussian matrix with rows 
$a_1,\dots,a_n$ then 
$$\frac{a_{1}}{\|a_{1}\|},\dots,\frac{a_{n}}{\|a_{n}\|},
\|a_{1}\|^{2},\dots,\|a_{n}\|^{2}$$ 
are independent random vectors and variables, the first $n$ of which
have uniform distribution on the sphere, $a_i/\|a_i\|\sim\UN
(\sphere^{m-1})$, and the last $n$ of which are $\chi^2_m$
distributed on $\R_+$. Recall that $\CN(A)$ is invariant 
under row scaling of $A$. The distribution of $\CN(A)$ under 
Gaussian input is therefore the same as under input matrices 
with i.i.d.\ $\UN(\sphere^{m-1})$ rows.

In subsequent sections we will therefore assume that 
\begin{equation*}
A_i:\Omega\rightarrow\sphere^{m-1}\qquad(i=1,\dots,n)
\end{equation*}
are i.i.d.\ random vectors defined on a probability space 
$(\Omega,{\mathcal F},\bfP)$ such that $A_i\sim\UN
(\sphere^{m-1})$, and that $A$ is the random $n\times m$ 
matrix 
\begin{equation*}
A=[A_1\dots A_n]^{\T}.
\end{equation*}
We say that $A$ is a {\em uniform random matrix}.
For convenience we shall assume that $m\leq n$, though the 
case $m<n$ is not difficult to derive from the results we will 
develop.

We adhere to the usual practice in probability theory and say 
that a property holds {\em almost surely} (or for almost all
$\omega\in\Omega$) if it holds with probability $1$. An event 
that occurs with probability zero is also called a null-set.

\begin{remark}\label{almost surely}
Using Proposition~\ref{support1} (ii) and Lemma~\ref{support3} 
it is trivial to show that if $A$ is a uniform random matrix then 
almost surely $A$ is not ill-posed. Likewise, the sets $S$ 
from \eqref{proof of support1} and $S'$ from \eqref{proof of support3} 
are bases of $\R^m$ almost surely. Finally, any $m$ i.i.d.\ uniformly 
drawn vectors from $\sphere^{m-1}$ are a basis of $\R^m$ almost surely.
Therefore, when $A$ is a uniform random matrix then for all 
$S\in{\mathcal P}_m$ the matrix $A_S$ is almost surely nonsingular 
and the random vectors $u_S$ and $\hat{u}_S$ are defined for almost 
all $\omega\in\Omega$.
\end{remark}
\hspace{1cm}

We denote the uniform probability measure on $\sphere^{m-1}$ by 
$\nu_{n-1}$. It is well known that the $m-1$-dimensional volume of 
$\Cap(p,\rho)$ equals the integral 
$I_{m-2}(\rho)=\int_0^r\sin^{m-2}x\;dx$ times the volume of the unit 
sphere in $\R^{m-1}$. Therefore, 
\begin{equation}\label{volume formula}
\bfP\left[A_i\in\Cap(p,\rho)\right]=
\nu_{m-1}(\Cap(p,\rho))=
\frac{I_{m-2}(\rho)}{I_{m-2}(\pi)}.
\end{equation}
It is trivial to check by induction that 
\begin{equation}\label{int_lemma}
I_m(\pi)\geq\frac{2}{\sqrt{m}}.
\end{equation}

\section{Upper Tail Bounds}\label{upper estimates}

The goal of this section is to derive an upper bound on the 
tail probability $\bfP[\CN(A)\geq t]$ when $A$ is a uniform random 
matrix.

Let $A$ be a random uniform matrix. For $S\in{\mathcal P}_m$ and 
$t\geq 1$ we consider the following events:
\begin{align*}
{\mathcal N}_t(S)&=\left\{\omega\in\Omega:\|u_{S}(\omega)\|\geq t
\right\},\\
{\mathcal B}^{\rm out}_t(S)&=\left\{\omega\in\Omega:A_i(\omega)
\notin\Cap\left(\hat{u}_{S}(\omega),
\arccos(1/t)\right)\forall i\notin S\right\},\\
{\mathcal B}^{\rm in}_{\pi/2}(S)&=\left\{\omega\in\Omega:A_i(\omega)
\in\Cap\left(\hat{u}_{S}(\omega),\pi/2\right)\forall i\notin
S\right\},\\
{\mathcal A}_t^{\rm if}&=\left\{\omega\in\Omega:\bigl(A(\omega)
\text{ infeasible }\bigr)
\wedge\bigl(\CN(A(\omega))\geq t\bigr)\right\},\\
{\mathcal A}_t^{\rm sf}&=\left\{\omega\in\Omega:\bigl(A(\omega)\text{ 
strictly feasible }\bigr)
\wedge\bigl(\CN(A(\omega))\geq t\bigr)\right\},\\
{\mathcal S}_t(S)&={\mathcal N}_t(S)\cap{\mathcal B}^{\rm out}_t,\\
{\mathcal S}_{t,\pi/2}(S)&={\mathcal N}_t(S)\cap{\mathcal
B}^{\rm in}_{\pi/2}.
\end{align*}

The following lemma is a key tool in our analysis.

\begin{lemma}\label{prere}
Let $A$ be a uniform random matrix and $t\in[1,\infty)$. 
Then $A_t^{\rm if}\setminus\bigl(\bigcup_{S\in{\mathcal P}_m}{\mathcal
S}_t(S)\bigr)$ and $A_t^{\rm sf}\setminus\bigl(\bigcup_{S\in{\mathcal
P}_m}{\mathcal S}_{t,\pi/2}(S)\bigr)$ are null-sets. 
\end{lemma}

\proof
If $A=A(\omega)$ is infeasible and $\CN(A)=1/\cos\theta(A)\geq t$ 
then 
\begin{equation}\label{theta ineq}
\arccos(1/t)\leq\theta(A)<\pi/2.
\end{equation} 
Let $\Cap\bigl(p,\theta(A)\bigr)$ be a LCP of $A$ and $S$ the 
corresponding blocking set. By Proposition~\ref{support1} and 
Remark~\ref{almost surely}, we have $|S|=m$ for almost all 
$\omega\in\Omega$. Lemma~\ref{support5} implies that $p(A)=\hat{u}_S$ 
for these $\omega$ and $\arccos(1/\|u_S\|)=\theta(A)$, which together 
with \eqref{theta ineq} shows that $\|u_S\|\geq t$ and 
\begin{equation*}
\Cap\left(\hat{u}_S,\arccos(1/t)\right)\subseteq
\Cap\left(\hat{u}_S,\arccos(1/\|u_S\|)\right)
\end{equation*}
does not contain any of the $A_i$ $(i\notin S)$. This shows that 
$A_t^{\rm if}\setminus\bigl(\bigcup_{S\in{\mathcal P}_m}{\mathcal
S}_t(S)\bigr)$ is a null-set.

If $A=A(\omega)$ is strictly feasible and $\CN(A)=1/|\cos
\theta(A)|\geq t$, then $\pi/2<\theta(A)\leq\arccos(-1/t)$ 
and hence, 
\begin{equation}\label{theta ineq 2}
\frac{\pi}{2}>\pi-\theta(A)\geq\arccos(1/t)
\end{equation}
Let $\Cap\bigl(-p(\omega),\pi-\theta(A))$ be the (unique) SCP of $A$. 
Lemma~\ref{support3} applied to this cap shows that there exists 
$\rho'\in[\pi-\theta(A),\pi/2)$ and $p'\in\sphere^{m-1}$ such that 
$\Cap(p',\rho')$ contains all rows of $A$. By Remark~\ref{almost
surely} $S=\bigl\{i:A_i(\omega)\in\partial\Cap(p',\rho')\bigr\}$ 
is of cardinality $m$ for almost all $\omega\in\Omega$. Lemma 
\ref{support5} implies that $p'=\hat{u}_S$ and
$\rho'=\arccos(1/\|u_S\|)$ for these $\omega$. Moreover, 
$\Cap(\hat{u}_S,\pi/2)\supset\Cap(\hat{u}_S,\rho')$ contains 
all rows of $A$ and in particular $\{A_i(\omega):i\notin S\}$. This 
shows that $A_t^{\rm sf}\setminus\bigl(\bigcup_{S\in{\mathcal
P}_m}{\mathcal S}_{t,\pi/2}(S)\bigr)$ is a null-set. 
\eproof
\hspace{1cm}

For shorthand notation we write $S^*=\{1,\ldots,m\}\in{\mathcal P}_m$,
${\mathcal S}^*_t={\mathcal S}_t(S^*)$ and ${\mathcal S}^*_{t,\pi/2}=
{\mathcal S}_{t,\pi/2}(S^*)$ in the sequel. Then, for $t\geq 1$, 
\begin{align}
\bfP\left[\CN(A)\geq t\right]&=\bfP\left[{\mathcal A}_t^{\rm if}\right]
+\bfP\left[{\mathcal A}_t^{\rm sf}\right]
+\bfP\left[\bigl(A\text{ ill-posed})\wedge\bigl(
\CN(A)\geq t\bigr)\right]\nonumber\\
&=\bfP\left[{\mathcal A}_t^{\rm if}\cap\bigcup_{S\in{\mathcal P}_m}
{\mathcal S}_t(S)\right]
+\bfP\left[{\mathcal A}_t^{\rm sf}\cap\bigcup_{S\in{\mathcal P}_m}
S_{t,\pi/2}(S)\right]+0\label{tri}\\
&\leq\sum_{S\in{\mathcal P}_m}\bfP\left[{\mathcal S}_t(S)\right]
+\sum_{S\in{\mathcal P}_m}\bfP\left[{\mathcal S}_{t,\pi/2}(S)\right]
\nonumber\\
&\leq\binom{n}{m}
\left(\bfP\left[{\mathcal S}^*_t\right]+\bfP\left[
{\mathcal S}^*_{t,\pi/2}\right]\right),\label{corol_cota}
\end{align}
where \eqref{tri} follows from Remark~\ref{almost surely} and 
Lemma~\ref{prere}.

Note that 
\begin{align}
\bfP\left[{\mathcal S}^*_t\right]&=
\bfP\left[{\mathcal B}^{\rm out}_t(S^*)\Bigl\|\bigr.
{\mathcal N}_t(S^*)\right]
P\left[{\mathcal N}_t(S^*)\right]\nonumber\\
=&\left(\int_{\sphere^{m-1}}\bfP\left[{\mathcal B}^{\rm out}_t
(S^*)\Bigl\|\bigr.{\mathcal N}_t(S^*),\hat{u}_{S^*}=x\right]
\nu_{m-1}(dx)\right)
\bfP\left[{\mathcal N}_t(S^*)\right]\label{uniform}\\
=&\left(\int_{\sphere^{m-1}}\left(1-\nu_{m-1}
\left(\Cap\left(x,\arccos(1/t)\right)\right)\right)^{n-m}
\nu_{m-1}(dx)\right)\cdot\bfP\left[
{\mathcal N}_t(S^*)\right]\label{independent}\\
=&\left(1-\frac{I_{m-2}\left(\arccos(1/t)\right)}
{I_{m-2}(\pi)}\right)^{n-m}\cdot
\bfP\left[{\mathcal N}_t(S^*)\right]\label{p8},
\end{align}
where \eqref{uniform} holds because 
$\DN\bigl(\hat{u}_{S^*}\bigl\|\bigr.{\mathcal N}_t(S^*)
\bigr)\sim \UN(\sphere^{m-1})$, 
\eqref{independent} holds because the random 
vectors $A_i$ $(i\notin S^*)$ are independent of $A_j$ $(j\in S^*)$, 
and \eqref{p8} follows from \eqref{volume formula}. Here, for a 
random variable $X$, $\DN(X)$ denotes the distribution of $X$, 
and $\UN(\sphere^{m-1})$ denotes the uniform distribution 
on $\sphere^{m-1}$. Likewise, a 
similar argument shows that 
\begin{equation}\label{p4}
\bfP\left[{\mathcal S}^*_{t,\pi/2}\right]=
2^{-(n-m)}\bfP\left[{\mathcal N}_t(S^*)\right].
\end{equation}

The task remains to determine bounds on $\bfP\left[{\mathcal N}_t(S^*)
\right]$. For $j\in S^*$ let $B_j$ be the unique unit vector in 
$\Span(\{A_i:i\neq j\})^{\perp}$ that complements $\{A_i:i\neq
j\}$ to a positively oriented basis of $\R^m$. Note that 
$A_j$ and $B_j$ are independent random vectors. 

\begin{lemma}\label{p1}
Let ${\mathcal C}_t:=\left\{\omega\in\Omega:\bigcup_{j\in S^*}
\left\{|B_j(\omega)
\cdot A_j(\omega)|\leq m/t\right\}\right\}$. Then 
${\mathcal C}_t\setminus{\mathcal N}_t(S^*)$ is a 
nullset for all $t\geq 1$. 
\end{lemma}

\proof
Almost surely $A_{S^*}$ is nonsingular and then  
\begin{equation}\label{ineq-frob}
\|u_{S^*}\|\leq\|A_{S^*}^{-1}\|\cdot\|e\|\leq\|A_{S^*}^{-1}\|_F\sqrt{m}.
\end{equation}
Together with \eqref{ineq-frob} the inequality 
$\|u_{S^*}\|\geq t$ implies that $\|A_{S^*}^{-1}\|_F\geq
t/\sqrt{m}$. Whenever this holds there exists an index $j\in S^*$ such that 
$\|A_{\cdot,j}^{-1}\|\geq t/m$, where $A_{\cdot,j}^{-1}$ denotes the 
$j$-th column of $A_{S^*}^{-1}$. Now the equation 
$A_{S^*}A_{S^*}^{-1}=I$ implies that $A_i\cdot A_{\cdot,j}^{-1}=
\delta_{ij}$ (the Kronecker symbol), which shows that 
$A_{\cdot,j}^{-1}/\|A_{\cdot,j}^{-1}=\pm B_j$ and 
$|B_j\cdot A_{j}|=1/\|A_{\cdot,j}^{-1}\|\leq m/t$. This proves 
the result. 
\eproof
\hspace{1cm}

\begin{lemma}\label{p2}
For all $m\geq 3$, $u\in\sphere^{m-1}$ and $(i=1,\dots,n)$ it is true 
that 
\begin{equation*}
\bfP\left[|A_i\cdot
u|\leq\frac{m}{t}\right]\leq\frac{m^{\frac{3}{2}}}{t}.
\end{equation*}
\end{lemma}

\proof 
Since the statement is trivially true when $m\geq t$, we may assume 
without loss of generality that $m<t$ and that $\arccos(m/t)$ is well 
defined. Therefore, 
\begin{equation*}
|A_i\cdot u|\leq m/t\Leftrightarrow
A_i\in\Cap\left(u,\arccos(-m/t)\right)\cap
\Cap\left(-u,\arccos(-m/t)\right).
\end{equation*}
Thus, 
\begin{align*}
\bfP\left[|A_i\cdot u|\leq m/t\right]
&=\left(I_{m-2}\left(\arccos(-m/t)\right)
-I_{m-2}\left(\arccos(m/t)\right)\right)\left(I_{m-2}(\pi)
\right)^{-1}\\
&\leq\dfrac{2\int_{\arccos(m/t)}^{\pi/2}
\sin x\;dx}{\int_0^\pi\sin^{m-2}x\;dx}
\leq\frac{m\sqrt{m-2}}{t}.
\end{align*}
where the last inequality follows from the fact that $m\geq 3$ and 
from equation \eqref{int_lemma}. 
\eproof
\hspace{1cm}

Lemma~\ref{p2} allows us to compute the bound we seek as follows:  
for $t\geq 1$ and $m\geq 3$ we have 
\begin{align}
\bfP\left[{\mathcal N}_t(S^*)\right] 
&\stackrel{\text{Lem~\ref{p1}}}{\leq}
\bfP\left[{\mathcal C}_t\right]+0\nonumber\\
&\leq\sum_{j\in S^*}
\bfP\left[|B_j\cdot A_{j}|\leq\frac{m}{t}\right]\nonumber\\
&=\sum_{j\in S^*}
\int_{\sphere^{m-1}}\bfP\left[|B_j\cdot A_{j}|\leq\frac{m}{t}
\Bigl\|\Bigr. B_j=x\right]\nu_{m-1}(dx)\label{use}\\
&\leq\frac{m^{\frac{5}{2}}}{t}.\label{p6}
\end{align}
where \eqref{use} uses the fact that $B_j$ is uniformly distributed 
on the sphere because the $A_i$ $(i\neq j)$ are, and where
\eqref{p6} follows from Lemma~\ref{p2} and the fact that $A_j$ is 
independent of $B_j$. Putting all the pieces together, we can now give 
an upper bound on the tail decay of $\CN(A)$. 

\begin{theorem}\label{tail_th}
For all $t\geq 1$, $m\geq 3$ and $n\geq m$ it is true that  
\begin{equation*}
\bfP\left[\CN(A)\geq t\right]
\leq\binom{n}{m}\cdot2m^{\frac{5}{2}}\cdot
\left(1-\frac{I_{m-2}\left(\arccos\left(\frac{1}{t}\right)
\right)}{I_{m-2}(\pi)}\right)^{n-m}\cdot\frac{1}{t}.
\end{equation*}
\end{theorem}

\proof
The claim follows immediately from equations \eqref{corol_cota}, 
\eqref{p8}, \eqref{p4} and \eqref{p6} together with 
the fact that $1-I_{m-2}(\arccos(1/t))/I_{m-2}(\pi)\geq 1/2$. 
\eproof
\hspace{1cm}

\section{Lower Tail Bounds}\label{lower estimates}

The goal of this section is to derive lower bounds on the decay 
rates of $\CN(A)$ in Theorem~\ref{lowbd}. In Section~\ref{exact decay} 
we will see that the combination of Theorems~\ref{tail_th} and 
\ref{lowbd} yields the exact asymptotic decay rates of $\log\CN(A)$.\\

Since $K\cap\sphere^{m-1}$ is a Borel set for all convex cones 
$K\subseteq\R^m$, the {\em angle space} of $K$ 
\begin{equation*}
\AS(K)=\nu_{m-1}(K\cap\sphere^{m-1})
\end{equation*}
is well-defined. It follows from the 
remarks of Section~\ref{probability model} that alternative equivalent 
definitions are provided by the relations 
\begin{equation*}
\AS(K)=\bfP\left[X\in K\cap\sphere^{m-1}\right]=
\bfP\left[Y\in K\right],
\end{equation*}
where $X\sim\UN(\sphere^{m-1})$ is a uniform random vector 
on the unit sphere and $Y$ is a multivariate normal random vector 
on $\R^m$ with covariance matrix $\sigma^2{\rm I}$ for any $\sigma^2>0$.

Let $A_i\sim\UN(\sphere^{m-1})$ $(i=1,\dots,k)$ be i.i.d.\ 
random vectors, where $k\leq m$, 
and let $\mathcal{LI}(A_1,\dots,A_k)\subset\Omega$ be the event that 
$\{A_1(\omega),\dots,A_{k}(\omega)\}$ is a linearly independent set 
of vectors. Then $\Omega\setminus
\mathcal{LI}(A_1,\dots,A_k)$ is a nullset, and for all 
$\omega\in\mathcal{LI}(A_1,\dots,A_k)$ there 
exists a unique orthogonal basis $\{E_1,\dots,E_k\}$ of 
$\Span\bigl(A_1,\dots,A_k\bigr)$ such that $E_i\cdot A_i>0$ and 
$\Span\bigl(E_1,\dots,E_i\bigr)=\Span\bigl(A_1,\dots,A_i\bigr)$ 
for $(i=1,\dots,k)$. In fact, the vectors $E_i$ are the column vectors 
of $Q$ in the thin $QR$ factorisation of the matrix $[A_1\dots A_k]$ 
(i.e., the $E_i$ are obtained by Gram-Schmidt orthogonalisation of the 
$A_i$). Let us consider the event 
\begin{equation}\label{C_m}
{\mathcal C}_{m,k}=\left\{\omega\in\mathcal{LI}(A_1,\dots,A_k):
\cone\bigl(A_1,\dots,A_k\bigr)\supseteq\cone\bigl(E_1,\dots,E_k\bigr)
\right\}.
\end{equation}
If $k=m$ we write ${\mathcal C}_{m}$ instead of ${\mathcal C}_{m,m}$. 
Note that in this case, 
$\AS\bigl(\cone(E_1,\dots,E_m)\bigr)=2^{-m}$. 
The following lemma shows thus that the angle space of the cone
generated by the $A_i$ is not too small with a quantifiable probability.

\begin{lemma}\label{IV}
Let $A_1,\dots,A_m$ be i.i.d.\ random vectors with
$A_i\sim\UN(\sphere^{m-1})$. Then 
\begin{equation}\label{equation I}
\bfP\left[{\mathcal C}_m\right]
\geq 2^{-\frac{2+m(m-1)}{2}},\qquad (m\geq 1). 
\end{equation}
\end{lemma}

\proof
We proceed by induction over $m$. For $m=1$ we have 
$\bfP\left[\cone(A_1)\supseteq\cone(E_1)\right]=1\geq 2^{-1}$, 
which shows that \eqref{equation I} holds true in the base case. 

Suppose \eqref{equation I} holds true for $m-1$ and let us show that 
it holds for $m$. For almost all $\omega\in\Omega$, 
\begin{equation}\label{frig}
A_m(\omega)\notin\Span\bigl(E_1(\omega),\dots,E_{m-1}(\omega)\bigr).
\end{equation} 
Let us thus assume that \eqref{frig} holds and let 
\begin{equation*}
\pi_{E_m}:\Span\bigl(E_1(\omega),\dots,E_m(\omega)\bigr)\rightarrow
\Span\bigl(E_1(\omega),\dots,E_{m-1}(\omega)\bigr)
\end{equation*}
denote the orthogonal projection along $E_m(\omega)$. Let 
\begin{equation*}
{\mathcal D}_m=\left\{\omega\in\Omega:
\frac{\pi_{E_m}(-A_m)}{\|\pi_{E_m}(-A_m)\|}\in\cone\left(E_1,\dots,E_{m-1}
\right)\right\}.
\end{equation*}
We claim that 
\begin{equation}\label{equation II}
{\mathcal D}_m\subseteq\left\{
E_m\in\cone\bigl(E_1,\dots,E_{m-1},A_m\bigr)\right\}.
\end{equation}
In fact, 
\begin{equation*}
\frac{\pi_{E_m}(-A_m)}{\|\pi_{E_m}(-A_m)\|}\in\cone\left(E_1,\dots,E_{m-1}
\right)\Leftrightarrow\pi_{E_m}(-A_m)\in\cone\left(E_1,\dots,E_{m-1}
\right),
\end{equation*}
except on a nullset, and $\pi_{E_m}(-A_m)\in\cone\left(E_1,\dots,E_{m-1}
\right)$ implies that there exist $\mu_i\geq 0$ $(i=1,\dots,m-1)$ 
and $\mu_m>0$ such that 
\begin{equation*}
A_m=\mu_m E_m+\pi_{E_m} A_m=\mu_m E_m-\sum_{i=1}^{m-1}
\mu_i E_i.
\end{equation*}
Hence, $E_m=\mu_m^{-1}\bigl(A_m+\sum_{i=1}^{m-1}\mu_i E_i\bigr)$, 
which proves \eqref{equation II}. Clearly \eqref{equation II} 
implies that 
\begin{equation}\label{equation III}
{\mathcal D}_m\cap{\mathcal C}_{m-1}\subseteq{\mathcal C}_m.
\end{equation}

Let $\G_{m-1,m}$ be the Grassmannian of $1$-codimensional linear 
subspaces of $\R^m$. $\G_{m-1,m}$ is a compact manifold with a
transitive group action defined by the orthogonal group $O_m$. 
The $G_{m-1,m}$-valued random variable 
\begin{equation}\label{G}
G:\omega\mapsto\Span(A_1(\omega),\dots,A_{m-1}(\omega)\bigr)\in
\G_{m-1,m}
\end{equation}
has uniform distribution $\nu_{G_{m-1,m}}$, that is, $\nu_{G_{m-1,m}}$ 
is the unique probability measure on $G_{m-1,m}$ that is invariant under 
the group action of $O_m$. In fact, this follows trivially from the 
spatial symmetry of the joint distribution of the $A_i$
$(i=1,\dots,m-1)$. It follows likewise from this symmetry that 
for all $g\in\G_{m-1,m}$ the random vectors 
\begin{equation*}
\frac{\pi_{E_m}(-A_m)}{\left\|\pi_{E_m}(-A_m)\right\|},A_1,\dots,
A_{m-1}
\end{equation*}
are independent random variables when conditioned on the event 
$\bigl\{\omega\in\Omega:G(\omega)=g\bigr\}$, with uniform conditional 
distributions 
\begin{equation}\label{eq7}
\DN\left(\frac{\pi_{E_m}(-A_m)}{\left\|\pi_{E_m}(-A_m)\right\|}
\Bigl\|\Bigr.G(\omega)=g\right),
\DN\left(A_i\Bigl\|\Bigr.G(\omega)=g\right)\sim\UN(\sphere^{m-1}\cap 
g)\equiv\UN(\sphere^{m-2}).
\end{equation}

These facts, \eqref{equation III} and the induction hypothesis finally 
imply 
\begin{align*}
\bfP&\left[{\mathcal C}_m\right]
\geq\bfP\left[{\mathcal D}_m\cap{\mathcal C}_{m-1}\right]\\
&\;=\int_{G_{m-1,m}}\Bigl(
\int_{y\in\otimes^{m-1}\sphere^{m-1}\cap g}
\bfP\Bigl[{\mathcal D}_m\Bigl\|\Bigr.(A_1,\dots,A_{m-1})=y,G=g\Bigr]
\otimes^{m-1}\nu_{m-2}(dy)\Bigr)\\
&\hspace{2cm}\cdot\bfP\left[{\mathcal C}_{m-1}
\Bigl\|\Bigr.G=g\right]\nu_{G_{m-1,m}}(dg)\\
&\;\geq\int_{G_{m-1,m}}\int_{y\in\otimes^{m-1}\sphere^{m-1}\cap g}
2^{-(m-1)}\otimes^{m-1}\nu_{m-2}(dy)\cdot\;2^{-\frac{2+(m-1)(m-2)}{2}}
\nu_{G_{m-1,m}}(dg)\\
&\;=2^{-\frac{2+m(m-1)}{2}}.
\end{align*}
\eproof

The next lemma shows that the angle space defined on a 1-codimensional 
hyperplane does not change too much under an orthogonal projection into 
a nearby hyperplane. 

\begin{lemma}\label{III}
Let $p_1,p_2\in\sphere^{m-1}$, let us denote the angle space 
defined on $p_i^{\perp}=\{x\in\R^m: p_i\cdot x=0\}$ by 
$\mathfrak{as}_i$ $(i=1,2)$, let $\pi_{p_2^{\perp}}$ be the 
orthogonal projection of $\R^m$ into $p_2^{\perp}$ along $p_2$ and 
let $\pi=\pi_{p_2^{\perp}}|_{p_1^{\perp}}$ be its restriction 
to $p_1^{\perp}$. Let finally $K$ be a convex cone in $p_1^{\perp}$. 
Then 
\begin{equation*}
\mathfrak{as}_2(\pi_2 K)\geq\mathfrak{as}_1(K)|p_1\cdot p_2|. 
\end{equation*}
\end{lemma}

\proof
If $p_1\cdot p_2=0$ then the bound is trivial. Therefore, w.l.o.g.\ 
$p_1\cdot p_2\neq 0$ and then $\pi$ is a vector space isomorphism 
between $p_1^{\perp}$ and $p_2^{\perp}$. Let $\{e_1,\dots,e_{m-2}\}$
be an orthonormal basis of $p_1^{\perp}\cap p_2^{\perp}$, let 
$e_j^{(i)}=e_j$ $(j=1,\dots,m-2)$ and let $e_{m-1}^{(i)}$ be chosen so
that $\{e_1^{(i)},\dots,e_{m-1}^{(i)}\}$ is an orthonormal basis of
$p_i^{\perp}$ for $(i=1,2)$. Then $|e_{m-1}^{(1)}\cdot e_{m-1}^{(2)}|
=|p_1\cdot p_2|$. 
Let us express vectors in $p_1^{\perp}$ in terms of coordinates 
$y\in\R^{m-1}$ defined by linear combinations $\sum_{j=1}^{m-1}
y_je_j^{(1)}$. Likewise, let $z$ be the coordinate system defined on 
$p_2^{\perp}$ by $\{e_1^{(2)},\dots,e_{m-1}^{(2)}\}$. 
Then $\pi$ expressed in terms of $y$--$z$ coordinates is the matrix 
\begin{equation*}
\pi=\begin{pmatrix}{\rm I}&0\\0&e_{m-1}^{(2)}\cdot e_{m-1}^{(1)}
\end{pmatrix}=\begin{pmatrix}{\rm I}&0\\0&\pm|p_1\cdot p_2|
\end{pmatrix}.
\end{equation*}
Now let $Z$ be a multivariate standard normal random vector on 
$p_2^{\perp}$, that is, $Z$ has the density function 
\begin{equation*}
f_{Z}(z)=(2\pi)^{-\frac{(m-1)}{2}}\exp\left(-\frac{\sum_{i=1}^{m-1}
z_i^2}{2}\right).
\end{equation*}
Then 
\begin{equation}\label{round1}
\mathfrak{as}_2(\pi K)=\bfP\left[Z\in\pi(K)\right]=\bfP\left[
Y\in K\right],
\end{equation}
where $Y=\pi^{-1} Z$ has density 
\begin{align*}
f_{Y}(y)&=f_{Z}(z(y))\left|\det\left(\frac{\partial z_i}{\partial y_j}
\right)\right|\\
&=(2\pi)^{-\frac{m-1}{2}}\exp\left(-\frac{1}{2}\left(
\sum_{i=1}^{m-2}y_i^2+y_{m-1}^2|p_1\cdot p_2|^2\right)
\right)|p_1\cdot p_2|\\
&\geq(2\pi)^{-\frac{m-1}{2}}\exp\left(-\frac{1}{2}
\sum_{j=1}^{m-1}y_i^2\right)|p_1\cdot p_2|,
\end{align*}
where the last inequality holds because $|p_1\cdot p_2|<1$. Therefore, 
\begin{equation*}
\bfP\left[Y\in K\right]\geq |p_1\cdot
p_2|\int_{K}(2\pi)^{-\frac{m-1}{2}}\exp\left(-\frac{\sum_{i=1}^{m-1}
y_i^2}{2}\right)dy=|p_1\cdot p_2|\mathfrak{as}_1(K).
\end{equation*}
Together with \eqref{round1} this proves the lemma. 
\eproof

The combination of Lemmas~\ref{support4},~\ref{IV} and~\ref{III} now 
allows us to derive lower bounds on the tail probabilities of $\CN(A)$.

\begin{theorem}\label{lowbd}
Let $A$ be a uniform random $n\times m$ matrix where $m\geq 2$ and 
$n\geq m$. Then there exists a constant $c(m)>0$ that depends only 
on $m$ such that for all $t\geq 1/\cos(\pi/4)$ it is true that 
\begin{equation*}
\bfP\left[\CN(A)\geq t\right]\geq c(m)\cdot
\left(\frac{I_{m-2}\left(\arccos\frac{1}{t}\right)}
{I_{m-2}(\pi)}\right)^{n-m}\cdot\frac{1}{t}.
\end{equation*}
\end{theorem}

\proof
It follows from \eqref{eq7}, the definition 
of $G$ in Lemma~\ref{IV} and the claim of the same result that 
\begin{align}
\bfP\left[{\mathcal C}_{m,k-1}\right]&=
\int_{G_{m-1,m}}\bfP\left[{\mathcal C}_{m,k-1}\|G=g\right]
\nu_{G_{m-1,m}}(dg)\nonumber\\
&\geq 2^{-\frac{2+(m-1)(m-2)}{2}}\int_{G_{m-1,m}}\nu_{G_{m-1,m}
}(dg)\nonumber\\
&=2^{-\frac{2+(m-1)(m-2)}{2}}.\label{EQUround1}
\end{align}
Since $\{A_1,\dots,A_m\}$ and $\{E_1,\dots,E_{m-1},A_m\}$ are linearly
independent sets for all $\omega\in{\mathcal C}_m$, it follows from 
Proposition~\ref{support1} (i) that the vectors 
$\{A_1,\dots,A_m\}$ define a unique SCP $\Cap\bigl(P_A(\omega),
R_A(\omega)\bigr)$, and likewise there exists a unique SCP 
$\Cap\bigl(P_E(\omega),R_E(\omega)\bigr)$ corresponding to the set of 
vectors 
$\{E_1,\dots,E_{m-1},A_m\}$. Moreover, for all $\omega\in{\mathcal C}_{m,k-1}$ 
we have $\Cap\bigl(P_E,R_E\bigr)\subseteq\Cap\bigl(P_A,R_A\bigr)$, and 
hence, 
\begin{equation}\label{EQUround2}
0\leq R_E(\omega)\leq R_A(\omega)\leq\frac{\pi}{2}\qquad
\forall\omega\in{\mathcal C}_{m,k-1}.
\end{equation}
Inequalities~\eqref{EQUround1} and~\eqref{EQUround2} 
imply that for all $t\geq 1$, 
\begin{align}
\bfP\left[R_A\geq\arccos 1/t\right]&\geq
\bfP\left[R_A\geq\arccos 1/t\|{\mathcal C}_{m,k-1}\right]
\bfP\left[{\mathcal C}_{m,k-1}\right]\nonumber\\
&\geq\bfP\left[R_E\geq
\arccos 1/t\|{\mathcal C}_{m,k-1}\right]\cdot 
2^{-\frac{2+(m-1)(m-2)}{2}}.\label{EQUround3}
\end{align}

Let $\{e_1,\dots,e_m\}$ be the canonical basis of $\R^m$. Then 
$\{e_1,\dots,e_{m-1},A_m\}$ is linearly independent almost surely, 
and it follows from Proposition~\ref{support1}~(i) that there
exists a unique SCP $\Cap(P,R)$ that corresponds to these vectors. 
Since $A_m$ is independent of the event ${\mathcal C}_{m,k-1}$, which is 
defined entirely in terms of $A_1,\dots,A_{m-1}$, and since the 
invariance of $\DN(A_i)$ under the action of the orthogonal group 
$O_m$ on $\sphere^{m-1}$ implies that $[E_1,\dots,E_{m-1}]$ is uniformly
distributed on the Stiefel manifold $V_{m-1}$ of $m\times(m-1)$
matrices with orthonormal columns, we have 
\begin{equation}\label{EQUround4}
\bfP\left[R_E\geq\arccos 1/t\bigl\|{\mathcal C}_{m,k-1}\right]
=\bfP\left[R\geq\arccos 1/t\right].
\end{equation}

We will consider the unit vectors
\begin{equation*}
\bar{e}=\frac{1}{\sqrt{m-1}}\sum_{i=1}^{m-1}e_i\quad\text{ and }
\quad p_\vartheta=\cos\vartheta\cdot e_m+\sin\vartheta\cdot \bar{e}
\quad\text{ for }\quad\vartheta\in[-\pi/2,\pi/2].
\end{equation*}

Then a random angle $\Theta:\Omega\rightarrow[-\pi/2,
\pi/2)$ is defined almost everywhere by the condition
$A_m(\omega)\in p_{\Theta(\omega)}^{\perp}+e_1$. In fact, 
\begin{equation*}
\Theta=\begin{cases}
\arctan\frac{\sqrt{m-1}A_m\cdot e_m}{1-\sum_{i=1}^{m-1}A_m\cdot e_i}
\qquad&\text{ if }\sum_{i=1}^{m-1}A_m\cdot e_i\neq 1,\\
-\frac{\pi}{2}&\text{ if }\sum_{i=1}^{m-1}A_m\cdot e_i=1\text{ and }
A_m\cdot e_m\neq 0,\\
\text{undefined}&\text{otherwise}.
\end{cases}
\end{equation*}
Note that $p_{\vartheta}^{\perp}+e_1=
p_{\vartheta}^{\perp}+e_i$ $(i=2,\dots,m-1)$ for all $\vartheta\in
[-\pi/2,\pi/2)$, that is, the definition of $\Theta$ is symmetric with
respect to the $e_i$.

It is easy to see that 
$\Theta$ has a continuous density function $f_{\Theta}>0$ such that 
$f_{\Theta}(-\vartheta)=f_{\Theta}(\vartheta)$ for all $\vartheta\in
(-\pi/2,\pi/2)$, and there exists a constant $c_{\Theta}>0$ such that 
$f_{\Theta}(\vartheta)\geq c_{\Theta}$ for all $\vartheta$ in the 
compact set $[-\pi/4,\pi/4]$. 

Note that one can parameterise the sphere $\sphere^{m-1}$ by 
$\sphere^{m-2}\times[-\pi/2,\pi/2)$ via 
\begin{equation*}
\left(\sum_{i=1}^{m-1}\lambda_i e_i,\vartheta\right)
\mapsto\left(\sum_{j=1}^m\mu_j e_j\right)\cdot\cos\vartheta
+\frac{1}{\sqrt{m-1}}p_{\vartheta},
\end{equation*}
for $\sum_{i=1}^{m-1}\lambda_i e_i\in\sphere^{m-2}\subset\Span(e_1,
\dots,e_{m-1})$, that is, $\sum_{i=1}^{m-1}\lambda_i^2=1$, and where 
\begin{equation*}
\begin{pmatrix}\mu_1\\ \vdots\\ \mu_m\end{pmatrix}=
\begin{pmatrix}
\frac{\cos\vartheta+(m-2)}{m-1}
&\frac{\cos(\vartheta)-1}{m-1}
&\dots
&\frac{\cos(\vartheta)-1}{m-1}
&\frac{\sin\vartheta}{\sqrt{m-1}}\\
\frac{\cos(\vartheta)-1}{m-1}
&\frac{\cos\vartheta+(m-2)}{m-1}
&
&\vdots
&\vdots\\
\vdots
&\vdots
&\ddots
&\vdots
&\vdots\\
\frac{\cos(\vartheta)-1}{m-1}
&\frac{\cos(\vartheta)-1}{m-1}
&\dots
&\frac{\cos\vartheta+(m-2)}{m-1}
&\frac{\sin\vartheta}{\sqrt{m-1}}\\
-\frac{\sin\vartheta}{\sqrt{m-1}}
&-\frac{\sin\vartheta}{\sqrt{m-1}}
&\dots
&-\frac{\sin\vartheta}{\sqrt{m-1}}
&\cos\vartheta
\end{pmatrix}
\begin{pmatrix}
\lambda_1\\ \vdots\\ \lambda_{m-1}\\0
\end{pmatrix}.
\end{equation*}
Note also that the matrix appearing in the display is orthogonal with 
last column corresponding to $p_{\vartheta}$. Thus, the chosen
parameterisation corresponds to tilting the unit sphere 
$\sphere^{m-2}\subset\Span(e_1,\dots,e_{m-1})$ by an angle $\vartheta$
about the affine hull $\aff(e_1,\dots,e_{m-1})$ and shrinking it 
by $\cos\vartheta$ to fit the radius of the sphere cut out of
$\sphere^{m-1}$ by the tilted plane. 

This parameterisation defines a conditional distribution 
$\DN(A_m\|\Theta=\vartheta)$ on $\sphere^{m-2}$ with continuous 
Radon-Nikodym derivative $f_{A_m\|\Theta}$ with respect to
$\nu_{m-2}$. Moreover, 
$f_{A_m\|\Theta}(x\|\vartheta)=0$ if and only if $x\in\aff(e_1,
\dots,e_{m-1})$. Therefore, there exists a constant $c_A>0$ such that 
\begin{equation}\label{EQUroundI}
f_{A_m\|\Theta}(x\|\vartheta)\geq c_A
\end{equation}
for all $(x,\vartheta)$ in the compact set 
$\left\{x\in\sphere^{m-1}:x\cdot\bar{e}\leq 0\right\}\times
\left[-\pi/4,\pi/4\right]$.

Now Lemma~\ref{support4} shows that for all $t\geq 1/\cos(\pi/4)$, 
\begin{align*}
\bfP\bigl[R\geq\arccos(1/t)\bigr]&=\bfP\bigl[
\bigl(\Theta\in
[\arccos 1/t-\pi/2,\pi/2-\arccos 1/t]\bigr)\\
&\hspace{3cm}\wedge\bigl(-\pi_{p_{\Theta}^{\perp}}A_m\in\cone(
\pi_{p_{\Theta}^{\perp}}e_1,\dots,\pi_{p_{\Theta}^{\perp}}e_{m-1})
\bigr))\bigr]\\
&=\int_{\arccos\frac{1}{t}-\frac{\pi}{2}}^{\frac{\pi}{2}
  -\arccos\frac{1}{t}}
\int_{\pi_{p_{\Theta}^{\perp}}(-\cone(e_1,\dots,e_{m-1}))}
f_{A_m\|\Theta}(x\|\vartheta)f_{\Theta}(\vartheta)\nu_{m-2}(dx)d\vartheta\\
&\kern-0.5cm\stackrel{\eqref{EQUroundI}, {\rm Lem~\ref{III}}}{\geq}
\int_{\arccos\frac{1}{t}-\frac{\pi}{2}}
^{\frac{\pi}{2}-\arccos\frac{1}{t}}c_A
\cdot 2^{-(m-1)}|\cos\vartheta|\cdot f_{\Theta}(\vartheta)
d\vartheta\\
&\geq c_A c_{\Theta}\cdot\left(2
\int_{0}^{\frac{\pi}{2}-\arccos\frac{1}{t}}\cos\vartheta
d\vartheta\right)\cdot 2^{-(m-1)}\\
&=c_A c_{\Theta}\cdot 2^{-(m-2)}\cdot\frac{1}{t}. 
\end{align*}
Therefore, \eqref{EQUround3} and \eqref{EQUround4} imply that 
\begin{equation}\label{rogue}
\bfP\left[R_A\geq\arccos 1/t\right]
\geq 2^{-\bigl(1+\frac{(m-1)(m-2)}{2}\bigr)}
\cdot c_A c_{\Theta}\cdot 2^{-(m-2)}\cdot\frac{1}{t}.
\end{equation}

Finally, if $A_{m+1},\dots,A_n\in\Cap(P_A,R_A)$ then $\Cap(P_A,R_A)$ 
is the SCP of $\{A_1,\dots,A_n\}$ and it follows from the remarks of 
Section~\ref{caps} that $\CN(A)=|\cos R_A|^{-1}$. Therefore, 
\begin{align*}
\bfP\left[\CN(A)\geq t\right]&\geq
\bfP\left[\left(A_{m+1},\dots,A_n\in\Cap(P_A,R_A)\right)
\wedge\left(R_A\geq\arccos 1/t\right)\right]\\
&\geq \left(\frac{I_{m-2}\left(\arccos\frac{1}{t}\right)}
{I_{m-2}(\pi)}\right)^{n-m}\cdot
2^{-\bigl(1+\frac{(m-1)(m-2)}{2}\bigr)}
\cdot c_A c_{\Theta}\cdot 2^{-(m-2)}\cdot\frac{1}{t}.
\end{align*}
Since $c_A$ and $c_{\Theta}$ depend only on $m$, this proves the 
claim of the theorem. 
\eproof

\section{Exact Tail Decay Rates}\label{exact decay}

The decay rates of $\bfP[\CN(A)\geq t]$ developed in
Sections~\ref{upper estimates} and~\ref{lower estimates} give an 
estimate on the rarity of a large backward error, high instability 
or long running times for some algorithms applied to a random 
linear feasibility problem drawn from Gaussian data.
Moreover, as mentioned in the introduction, the best upper bounds 
on the running time of modern linear programming or linear feasibility 
solvers applied to real input data are polynomial in the problem 
dimension and $\log\CN(A)$. We are therefore also interested in 
estimates on probability tails 
\begin{equation}\label{logdecay}
\bfP[\log\CN(A)\geq t]
\end{equation}
for $t\gg 1$. 

Theorem~\ref{tail_th} implies that 
\begin{equation}\label{round A}
\bfP[\log\CN(A)\geq t]=\bfP[\CN(A)\geq\e^{t}]\leq
\binom{n}{m}2
m^{\frac{5}{2}}\left(1-\frac{I_{m-2}(\arccos(\e^{-t}))}
{I_{m-2}(\pi)}\right)^{n-m}\cdot\e^{-t}.
\end{equation}
On the other hand, Theorem~\ref{lowbd} shows 
\begin{equation*}
\bfP[\log\CN(A)\geq t]\geq
c(m)\left(\frac{I_{m-2}(\arccos(\e^{-t}))}{I_{m-2}(\pi)}
\right)^{n-m}\cdot\e^{-t}. 
\end{equation*}
Since $I_{m-2}(\arccos\e^{-t})/I_{m-2}(\pi)$ increases
monotonically to $1/2$ for $t\rightarrow\infty$, these formulas
show that the exponential decay rate of \eqref{logdecay} is exactly 
$-1$.

\begin{corollary}\label{exact decay rates}
If $A$ is a random uniform $n\times m$ matrix then 
\begin{equation*}
\lim_{t\rightarrow\infty}\frac{\log\bfP[\log\CN(A)\geq t]}
{t}=-1.
\end{equation*}
\end{corollary}

\proof
The proof is immediate from the arguments above. 
\eproof

Thus, although the multiplicative constant in \eqref{round A} is too
large, the formula captures the correct qualitative behaviour of the 
tails of $\log\CN(A)$ and the best possible upper bound on
\eqref{logdecay} must be of the form 
\begin{equation}\label{round B}
\bfP[\log\CN((A)\geq t]\leq c(m,n)\cdot\e^{-t}
\end{equation} 
for some constant $c(m,n)$ that depends on $m$ and $n$.

The exponential decay of $\bfP[\log\CN((A)\geq t]$ shows that 
the linear feasibility problem, and by extension linear programming, 
is ``empirically strongly polynomial''. See Section \ref{discussion} 
for further comments on this important point. 

\section{Moment Estimates}\label{moment estimates}

The probabilistic analysis of linear programming is primarily
concerned with the average running time of LP algorithms on random
input data. Because complexity bounds for interior-point methods 
are polynomial in $\log \CN(A)$ (see the introduction), upper bounds 
on the expectation, the variance and higher moments of the running time 
are easily derived from upper bounds on the corresponding moments of 
$\log\CN(A)$. 

Since $\bfE[X]=\int_{0}^{\infty}\bfP[X>x]dx$ for any random variable
$X$ that takes only nonnegative values, the estimate \eqref{round A} 
can be used to derive upper bounds on all moments of $\log\CN(A)$. 
Indeed, \eqref{round B} shows that for all $\gamma>0$,  
\begin{equation*}
\int_{0}^{\infty}\bfP\left[\left(\log\CN(A)\right)^{\gamma}\geq t\right]
d t\leq\int_{0}^{\infty}c(n,m)\e^{-t^{\frac{1}{\gamma}}}d t=
c(n,m)\Gamma(\gamma+1)<\infty,
\end{equation*}
that is, all moments of $\CN(A)$ are finite. To turn this into a 
quantitative estimate, we consider the function 
$\varphi:\R_+\to(\frac12,1]$ defined as follows:
$$
  \varphi(t)=1-\frac{I_{m-2}(\arccos(\e^{-t}))}{I_{m-2}(\pi)}.
$$
Note that $\varphi$ is continuous, and strictly decreasing with 
$\varphi(0)=1$ and 
$\displaystyle\lim_{t\to\infty}\varphi(t)=\frac12$. Let us define 
$$
f(m,n)=\varphi^{-1}\left((1/2)^{1/\sqrt{n}}\right). 
$$
Since $f(m,n)>0$, the following result follows.

\begin{corollary}\label{estimates}
Let $A$ be a uniform random $n\times m$ matrix with $n\geq m\geq 3$. 
Then for all $\gamma\in\R_+$ the $\gamma$-th moment of 
$\log\CN(A)$ is bounded by 
\begin{equation*}
\bfE\left[\left(\log\CN(A)\right)^{\gamma}\right]\leq
f(m,n)^{\gamma}+\binom{n}{m} 2 m^{\frac{5}{2}} 2^{-\frac{n-m}{\sqrt{n}}}
\Gamma(\gamma+1).
\end{equation*}
\end{corollary}

\proof
Using \eqref{round A}, we find
\begin{align*}
\bfE\left[\left(\log\CN(A)\right)^{\gamma}\right]&=
\int_{0}^{\infty}\bfP\left[\left(\log\CN(A)\right)^{\gamma}>t\right]d t\\
&\leq f(m,n)^{\gamma}+\int_{f(m,n)^{\gamma}}
^\infty\binom{n}{m}2 m^{\frac{5}{2}}
\left(1-\frac{I_{m-2}\bigl(\arccos(\e^{-f(m,n)})\bigr)}{I_{m-2}(\pi)}
\right)^{n-m}\cdot\e^{-t^{\frac{1}{\gamma}}}d t\\
&\leq f(m,n)^{\gamma}+\binom{n}{m}2m^{\frac{5}{2}}2^{-\frac{n-m}{\sqrt{n}}}
\int_{f(m,n)^{\gamma}}^\infty \e^{-t^{\frac{1}{\gamma}}}d t\\
&\leq f(m,n)^{\gamma}+\binom{n}{m}2m^{\frac{5}{2}}2^{-\frac{n-m}{\sqrt{n}}}
\Gamma(\gamma+1). 
\end{align*}
\eproof
\hspace{1cm}

In Section~\ref{limit theorems} we will see that the bounds of 
Corollary~\ref{estimates} are particularly useful for understanding
the behaviour of $\CN(A)$ when $n\gg m$. Note however that these 
bounds grow exponentially in $m$. One of the major objectives of the 
probabilistic analysis of linear programming is to show that the 
expected running 
times of particular families of algorithms are bounded by a polynomial 
of the dimension of the input data. Such results are often interpreted 
in the light of ``average strongly polynomiality'' of linear
programming.

Does the exponential growth of the estimates from Corollary~\ref{estimates} 
mean that Theorem~\ref{tail_th} fails to lead to ``average strong 
polynomiality'' results when used to bound the complexity of
interior-point algorithms for example? 
Not in the least! The exponential growth of 
the estimates from Corollary~\ref{estimates} is purely a consequence of 
our definition of the cut-off point $f(m,n)$, which we chose so as to 
converge to zero as $n$ tends to infinity to fit the purposes of 
the limit theorems of Section~\ref{limit theorems}. Giving up on this 
condition one can easily derive linear bounds:

\begin{lemma}\label{polylemma}
Let $\bigl\{X_{m,n}: (m,n)\in\N\times\N\bigr\}$ be a set of random
variables and $\gamma\geq 1$ a real number. Furthermore, let 
$p(m,n)$ and $t(m,n)$ be functions of $m$ and $n$ such that 
\begin{equation*}
\bfP\left[X_{m,n}>t\right]\leq\e^{p(m,n)-t^{\frac{1}{\gamma}}}
\qquad\forall\,t\geq t(m,n).
\end{equation*}
Then 
\begin{equation*}
\bfE\left[X_{m,n}\right]\leq \max\bigl(p(m,n)^{\gamma},t(m,n)\bigr)+
\Gamma(\gamma+1)2^{\gamma-1}.
\end{equation*}
\end{lemma}

\proof
\begin{align}
\bfE\left[X\right]&=\int_{0}^{\infty}\bfP\left[X>t\right]dt
\nonumber\\
&\leq\int_{0}^{\max(p(m,n)^{\gamma},t(m,n))}1 dt
+\int_{\max(p(m,n)^{\gamma},t(m,n))}^{\infty}\exp\left\{p(m,n)
-t^{\frac{1}{\gamma}}\right\}dt\nonumber\\
&\leq\max\bigl(p(m,n)^{\gamma},t(m,n)\bigr)+\int_{0}^{\infty}
\exp\left\{-t^{\frac{1}{\gamma}}2^{1-\frac{1}{\gamma}}\right\}dt
\label{EXPLAIN}\\
&=\max\bigl(p(m,n)^{\gamma},t(m,n)\bigr)+
\Gamma(\gamma+1)2^{\gamma-1},\nonumber
\end{align}
where \eqref{EXPLAIN} holds true because we claim that 
\begin{equation*}
p(m,n)-t^{\frac{1}{\gamma}}2^{1-\frac{1}{\gamma}}\leq
-\left(t-\max\bigl(p(m,n)^{\gamma},t(m,n)\bigr)
\right)^{\frac{1}{\gamma}}
\end{equation*}
for all $t\geq\max\bigl(p(m,n)^{\gamma},t(m,n)\bigr)$.
In fact, $x\mapsto x^{\frac{1}{\gamma}}$ is a concave function, 
since $\gamma\geq 1$. Therefore, 
\begin{align*}
\frac{1}{2}\left(p(m,n)^{\gamma}\right)^{\frac{1}{\gamma}}
&+\frac{1}{2}\left(t-\max\bigl(p(m,n)^{\gamma},t(m,n)\bigr)
\right)^{\frac{1}{\gamma}}\\
&\leq 2^{-\frac{1}{\gamma}}\left(p(m,n)^{\gamma}+t-
\max\bigl(p(m,n)^{\gamma},t(m,n)\bigr)\right)^{\frac{1}{\gamma}}\\
&\leq  2^{-\frac{1}{\gamma}}t^{\frac{1}{\gamma}},
\end{align*}
which shows that 
\begin{equation*}
p(m,n)+\left(t-\max\bigl(p(m,n)^{\gamma},t(m,n)\bigr)\right)
^{\frac{1}{\gamma}}\leq 2^{1-\frac{1}{\gamma}}t^{\frac{1}{\gamma}}
\end{equation*}
and proves our claim. 
\eproof
\hspace{1cm}

\begin{corollary}\label{polycor}
Let $A$ be a uniform random $n\times m$ matrix where $n\geq m\geq 3$, 
and let $\gamma\geq 1$ be a real number. Then 
\begin{equation*}
\bfE\left[\left(\log\CN(A)\right)^{\gamma}\right]\leq 
\left(m\log n+\frac{5}{2}\log m+\log
2\right)^{\gamma}+\Gamma(\gamma+1)
2^{\gamma-1}.
\end{equation*}
In particular, 
\begin{align*}
\bfE\left[\log\CN(A)\right]&\leq m\log n+\frac{5}{2}\log m+\log 2+1,\\
\VAR\left(\log\CN(A)\right)&\leq\left(m\log n+\frac{5}{2}\log m+\log
2\right)^2+4.
\end{align*}
\end{corollary}

\proof
Equation \eqref{round A} shows that for all $t\geq t(m,n)=1$, 
\begin{align*}
\bfP\left[\left(\log\CN(A)\right)^{\gamma}>t\right]&=
\bfP\left[\log\CN(A)>t^{\frac{1}{\gamma}}\right]\\
&\leq\binom{n}{m}2 m^{\frac{5}{2}}
\left(1-\frac{I_{m-2}\bigl(\arccos(\exp\{-t^{\frac{1}{\gamma}}\})
\bigr)}{I_{m-2}(\pi)}\right)^{n-m}\cdot\e^{-t^{\frac{1}{\gamma}}}\\
&\leq n^m 2 m^{\frac{5}{2}}\e^{-t^{\frac{1}{\gamma}}}\\
&\leq\e^{m\log n+\frac{5}{2}\log m+\log 2-t^{\frac{1}{\gamma}}}.
\end{align*}
The first claim now follows from Lemma~\ref{polylemma}. Finally, since 
\begin{equation*}
\VAR(\log\CN(A))=\bfE[(\log\CN(A))^2]-\bfE[\log\CN(A)]^2\leq
\bfE[(\log\CN(A))^2],
\end{equation*}
the last two claims are special cases of the first claim.
\eproof
\hspace{1cm}

Corollary~\ref{polycor} recovers the main result 
in~\cite{ChC01}. However, it still does not fully exhaust the potential 
power of Theorem \ref{tail_th} and Lemma~\ref{polylemma}: indeed, in 
Section~\ref{mgg1} we will further strengthen Corollary~\ref{polycor} 
and show that $\bfE[\log\CN(A)]$ is asymptotically bounded by $m\log 2$ 
for arbitrary $n\geq m$, and by $\Oh(m^{\gamma})$ for any $\gamma>0$ 
when $n\geq 5m$.

\begin{remark}\label{suspicion}
Let us briefly remark here that the main reason for the appearance of 
the binomial term $\binom{n}{m}$ in the bound of Theorem~\ref{tail_th} 
is a lack of proper understanding of $\DN(A_i\|P=p)$, where $P$ is the 
centre of an LCP of $A=[A_1,\dots,A_n]^{\T}$. We suspect that 
knowledge of this conditional distribution would make it possible to 
replace $\binom{n}{m}$ by a polynomial term in $m$ and  $n$. If this 
hunch were true then the bound of Corollary~\ref{polycor} would of 
course become logarithmic in $m$ and $n$. 
\end{remark}

Let us finally investigate the moments of $\CN(A)$ itself, which is 
interesting in its own right for reasons mentioned in the introduction. 

\begin{corollary}\label{CORI}
Let $A$ be a uniform random $n\times m$ matrix where $n\geq m\geq 3$. 
Then 
\begin{equation*}
\bfE[\CN(A)^{\gamma}]\begin{cases}
=+\infty\qquad&\text{ if }\gamma\geq 1,\\
\leq 1+\binom{n}{m}2 m^{\frac{5}{2}}\frac{\gamma}{1-\gamma}\quad&
\text{ if }\gamma\in(0,1).
\end{cases}
\end{equation*}
\end{corollary}

\proof
Theorem~\ref{lowbd} shows that for all $t\geq(\cos(\pi/4))^{-\gamma}$, 
\begin{equation*}
\bfP\left[\CN(A)^{\gamma}>t\right]\geq
c(m)\left(\frac{I_{m-2}\bigl(\arccos t^{-\frac{1}{\gamma}}\bigr)}
{I_{m-2}(\pi)}\right)^{n-m}t^{-\frac{1}{\gamma}}.
\end{equation*}
Therefore, 
\begin{align*}
\bfE\left[\CN(A)^\gamma\right]&=\int_{0}^{\infty}
\bfP\left[\CN(A)^{\gamma}>t\right]dt\\
&\geq\int_{\bigl(\cos(\frac{\pi}{4})\bigr)^{-\gamma}}^{\infty}
c(m)2^{m-n} t^{-\frac{1}{\gamma}}dt=\infty\qquad\forall \gamma\geq 1.
\end{align*}
On the other hand, using Theorem
\ref{tail_th}, we find that for $\gamma\in(0,1)$, 
\begin{align*}
\bfE\left[\CN(A)^{\gamma}\right]&=
1+\int_{1}^{\infty}\bfP\left[\CN(A)^{\gamma}>t\right]dt\\
&\leq 1+\binom{n}{m}2m^{\frac{5}{2}}\int_{1}^{\infty}t^{-\frac{1}
{\gamma}}dt\\
&=1+\binom{n}{m}2m^{\frac{5}{2}}\frac{\gamma}{1-\gamma}.
\end{align*}
\eproof

Note that the moment bound of Corollary~\ref{CORI} grows exponentially
in $n$ for $\gamma<1$. This bound does not reflect the correct
limiting behaviour, since we will show in Corollary~\ref{corollary 6} 
below that $\lim_{n\rightarrow\infty}\bfE\bigl[
\CN(A)^\gamma\bigr]=1$ occurs for $\gamma<1$ and $m$ fixed.

\section{Limit Theorems for $n\gg m$}\label{limit theorems}

In this section we investigate the behaviour of $\CN(A)$ and 
$\log\CN(A)$ in the situation where $n\gg m$. For example, we 
will show that for fixed $m$, 
\begin{equation}\label{rround 1}
\CN(A)\stackrel{n\rightarrow\infty}{\longrightarrow}1
\end{equation}
with probability $1$, and 
\begin{equation}\label{rround 2}
\bfE\left[\log\CN(A)\right]\stackrel{n\rightarrow\infty}
{\longrightarrow}0.
\end{equation}
Intuitively it is clear that this behaviour should be observed:
whenever the system $Ax\leq 0, x\neq 0$ contains a very large numbers
of random constraints, the system should be infeasible and this 
infeasibility should be easy to detect algorithmically. Our results
confirm this intuition.\\

In the results below $m\geq 3$ is a fixed dimension, 
$(A_i)_{\N}$ denotes a sequence of i.i.d.\ random 
vectors with uniform distribution on the sphere, 
$\DN(A_i)\sim\UN(\sphere^{m-1})$, and $(A^{[n]})_{\N}$ is the sequence 
of random matrices $A^{[n]}=[A_1,\dots,A_n]^{\T}$. 

\begin{theorem}\label{pr1}
Let $(A_i)_{\N}$ and $(A^{[n]})_{\N}$ be as defined above. Then 
\begin{equation*}
\bfP\left[\lim_{n\rightarrow\infty}\CN(A^{[n]})=1\right]=1.
\end{equation*}
\end{theorem}

\proof
For all $\omega\in\Omega$ and $n\in\N$ let
$\Cap(P_n(\omega),R_n(\omega))$ be a LCP of $A^{[n]}$. By virtue of 
Proposition~\ref{geo3} it suffices to prove that 
\begin{equation}\label{round Q}
\bfP\left[\lim_{n\rightarrow\infty}R_n=0\right]=1.
\end{equation}
Let $\rho\in(0,\pi/2)$ be a fixed radius. Since $\sphere^{m-1}$ is 
compact, there exists a finite set of vectors $\{p_1,\dots,p_k\}\subset
\sphere^{m-1}$ such that $\bigcup_{i=1}^k\Cap(p_i,\rho/2)=\sphere^{m-1}$. 
By 
\begin{equation*}
\mathcal{CE}_{i,n}=\left\{\omega\in\Omega:A_1,\dots,A_n\notin
\Cap(p_i,\rho/2)\right\}
\end{equation*}
let us denote the event that the $i$-th cap does not contain any of 
the $n$ first vectors of $(A_i)_{\N}$. Then 
\begin{equation*}
\bfP\left[\mathcal{CE}_{i,n}\right]=\left(\frac{I_{m-2}(\pi-\rho/2)}
{I_{m-2}(\pi)}\right)^n,
\end{equation*}
and hence, 
\begin{equation}\label{round R}
\bfP\left[\bigcup_{i=1}^k\mathcal{CE}_{i,n}\right]\leq
\sum_{i=1}^k\bfP[\mathcal{CE}_{i,n}]=
k\cdot\left(\frac{I_{m-2}(\pi-\rho/2)}{I_{m-2}(\pi)}\right)^n.
\end{equation}

We now claim that 
\begin{equation}\label{round S}
\left\{\omega\in\Omega:R_n\geq\rho\right\}\subseteq
\bigcup_{i=1}^k\mathcal{CE}_{i,n}.
\end{equation}
In fact, if $\omega\in\bigl(\bigcup_{i=1}^k\mathcal{CE}_{i,n}\bigr)^c$, 
the complement of $\bigcup_{i=1}^k\mathcal{CE}_{i,n}$,  
then there exist indices $i\in\{1,\dots,k\}$ and $j\in\{1,\dots,n\}$
such that $P_n(\omega)\in\Cap(p_i,\rho/2)$ and
$A_j(\omega)\in\Cap(p_i,\rho/2)$. Using the triangular inequality on
the sphere we find $R_n\leq\arccos\bigl(P_n(\omega)\cdot
A_j(\omega)\bigr)<2\cdot\rho/2$. This shows 
\begin{equation*}
\left\{\omega\in\Omega:R_n<\rho\right\}\supseteq\left(
\bigcup_{i=1}^k\mathcal{CE}_{i,n}\right)^c,
\end{equation*} 
which is equivalent to our claim. 

Now \eqref{round R} and \eqref{round S} show 
\begin{equation*}
\bfP\left[R_n\geq\rho\right]\leq k\cdot\left(\frac{I_{m-2}
(\pi-\rho/2)}{I_{m-2}(\pi)}\right)^n\stackrel{n\rightarrow\infty}
{\longrightarrow}0.
\end{equation*}
Finally, since $n_1\leq n_2$ implies $R_{n_1}\geq R_{n_2}$, we have 
\begin{equation*}
0\leq\bfP\left[\lim_{n\rightarrow\infty} R_n>\rho\right]\leq
\lim_{n\rightarrow\infty}\bfP\left[R_n>\rho\right]=0.
\end{equation*}
Since this is true for all $\rho\in(0,\pi/2)$, \eqref{round Q} follows.
\eproof

\begin{corollary}\label{corollary 4}
Let $(A^{[n]})_{\N}$ be as above. Then
\begin{itemize}
\item[i)\;]
$\CN(A^{[n]})\stackrel{n\rightarrow\infty}{\longrightarrow_P} 1$,
\item[ii)\;]
$\CN(A^{[n]})\stackrel{n\rightarrow\infty}{\Longrightarrow} 1$,
\item[iii)\;] 
$\displaystyle\bfP\left[\lim_{n\rightarrow\infty}
\log\CN(A^{[n]})=0\right]=1$,
\item[iv)\;]
$\log\CN(A^{[n]})\stackrel{n\rightarrow\infty}{\longrightarrow_P} 0$,
\item[v)\;]
$\log\CN(A^{[n]})\stackrel{n\rightarrow\infty}{\Longrightarrow} 0$,
\end{itemize}
where $\longrightarrow_P$ denotes convergence in probability and 
$\Rightarrow$ denotes weak convergence. 
\end{corollary} 

\proof
These are all standard consequences of Theorem~\ref{pr1}, see for
example Theorem~25.2 in~\cite{Billingsley}. 
\eproof
\hspace{1cm}

Using \eqref{round A} and Corollary~\ref{corollary 4} one can analyse 
the asymptotic behaviour of $\bfE[\CN(A)]$ and $\bfE[\log\CN(A)]$, 
see Corollary~\ref{corollary 6} below. In the case of $\log\CN(A)$ 
for example, one can show that 
$\lim_{n\rightarrow\infty}\bfE[\log\CN(A^{[n]})]=0$, using Skorhohod's 
theorem. Remarkably, the estimates of Corollary~\ref{estimates} are strong 
enough to yield this result directly, without resort to 
Theorem~\ref{pr1}.

\begin{corollary}\label{limits}
Let $m$ be fixed and $(A_i)_{\N}$ and $(A^{[n]})_{\N}$ defined 
as above. Then 
\begin{equation*}
\lim_{n\rightarrow\infty}\bfE\left[(\log\CN(A^{[n]}))^{\gamma}
\right]=0,\qquad\forall \gamma>0.
\end{equation*}
In particular, 
\begin{itemize}
\item[i)\;]
$\displaystyle\lim_{n\rightarrow\infty}\bfE\bigl[\log\CN(A^{[n]})\bigr]=0$ 
and
\item[ii)\;]
$\displaystyle\lim_{n\rightarrow\infty}\VAR\left(\log\CN(A^{[n]})\right)=0$.
\end{itemize}
\end{corollary}

\proof
Since $2^{-(n-m)/\sqrt{n}}$ is exponentially decreasing in $n$ and 
\begin{equation*}
\binom{n}{m}\leq n^m
\end{equation*}
increases only polynomially in $n$, the second term in the estimate of 
Corollary~\ref{estimates} tends to zero as $n$ tends to infinity. 
In addition, $f(m,n)$ tends to zero as $n$ tends to infinity by 
definition of $f(m,n)$. This proves the displayed formula and 
as a particular case part {\em i)}. Part {\em ii)} follows from 
the display and the fact that $\VAR(\log\CN(A))\leq\bfE[(\log\CN(A))^2]$. 
\eproof
\hspace{1cm}

Let us now analyse the asymptotic behaviour of $\bfE\left[\CN(A)\right]$. 
Recall from Corollary~\ref{CORI} that $\bfE[\CN(A)^{\gamma}]$ is
finite if and only if $\gamma<1$. 

\begin{corollary}\label{corollary 6}
Let $m$ be fixed and $(A_i)_{\N}$ and $(A^{[n]})_{\N}$ defined 
as above. Then
\begin{equation*}
\lim_{n\rightarrow\infty}\bfE\bigl[(\CN(A^{[n]}))^{\gamma}\bigr]=1,
\qquad\forall\gamma\in[0,1).
\end{equation*}
\end{corollary}

\proof
The result is trivial for $\gamma=0$. Let us therefore assume that 
$\gamma\in(0,1)$. 
Let $t_{0}\in\R_+$ be large enough so that 
\begin{equation*}
1-\frac{I_{m-2}\left(\arccos t^{-\frac{1}{\gamma}}
\right)}{I_{m-2}(\pi)}\leq
\frac{3}{4},\qquad\forall t\geq t_{0}.
\end{equation*}
Since $\CN(A^{[n]})\rightarrow_{P}1$ by Corollary~\ref{corollary 4}, 
for all $\epsilon>0$ there exists a number $n_{\epsilon}\in\N$ such 
that 
\begin{equation}\label{in probability}
\bfP\left[\CN(A^{[n]})>(1+\epsilon)^{\frac{1}{\gamma}}\right]
<\epsilon,\qquad\forall n\geq n_{\epsilon}.
\end{equation}
Therefore, for $n\geq n_{\epsilon}$, 
\begin{align*}
\bfE&\left[\left(\CN(A^{[n]})\right)^{\gamma}\right]=\int_{0}^{\infty}
\bfP\left[\left(\CN(A^{[n]})\right)^{\gamma}>t\right]dt\\
&\leq\int_{0}^{(1+\epsilon)^{\frac{1}{\gamma}}}1dt+
\int_{(1+\epsilon)^{\frac{1}{\gamma}}}^{t_0}
\bfP\left[\CN(A^{[n]})>(1+\epsilon)^{\frac{1}{\gamma}}\right]dt
+\int_{t_0}^{\infty}\bfP\left[\CN(A^{[n]})>t^{\frac{1}{\gamma}}\right]dt\\
&\stackrel{\text{Thm}\ref{tail_th},\eqref{in probability}}{\leq} 
(1+\epsilon)^{\frac{1}{\gamma}}+\epsilon\left(t_0-(1+\epsilon)^{\frac
{1}{\gamma}}\right)+
\binom{n}{m}2m^{\frac{5}{2}}\left(\frac{3}{4}\right)^{n-m}
\int_{t_0}^\infty t^{-\frac{1}{\gamma}}dt\\
&=(1+\epsilon)^{\frac{1}{\gamma}}+\epsilon
\left(t_0-(1+\epsilon)^{\frac{1}{\gamma}}\right)+
2n^m m^{\frac{5}{2}}\left(\frac{3}{4}\right)^{n-m}\frac{t_0^{1
-\frac{1}{\gamma}}}{\frac{1}{\gamma}-1}.
\end{align*}
Taking limits as $n\rightarrow\infty$ and observing that $\epsilon>0$ 
was arbitrary, the claim follows. 
\eproof
\hspace{1cm}

\section{Limit Theorems for $m\gg 1$}\label{mgg1}

In this section we investigate the behaviour of $\log\CN(A)$ in the 
situation where $m\gg 1$. We will see that 
$\lim\sup_{m\rightarrow\infty}\bfE[\log\CN(A)]/m
\leq\log 2$, and we point out why we suspect that the correct value 
of this limit is zero.

Let $(\rho_m)_{\N}\subset(0,\pi/2]$ be a sequence such that 
$\displaystyle\lim_{m\rightarrow\infty}\rho_m=\pi/2$, and 
let $(A^{[m]})_{\N}$ be a sequence of random vectors such that 
$A^{[m]}\sim\UN(\sphere^{m-1})$. It can be shown that if $\rho_m$ 
converges to $\pi/2$ at an algebraic rate as a function of $m$, then 
\begin{equation*}
\lim_{m\rightarrow\infty}\bfP\left[A^{[m]}\in\Cap(p,\rho_m)\right]=0.
\end{equation*}
This effect is a special case 
of the so-called {\em concentration of measure phenomenon}, see e.g.\ 
\cite{Ledoux}. The phenomenon is 
remarkable, because it implies that after fixing an equator by choosing an 
arbitrary grand circle on a high dimensional sphere, one will 
observe that a counter-intuitively high proportion of uniformly drawn 
sample points from that sphere lie in a very narrow neighbourhood 
around that equator. This phenomenon affects the analysis of the 
distribution tails of $\CN(A)$ for large $m$.

However, for the purposes of this analysis it suffices to know that 
for any fixed real exponent $\gamma>0$,  
\begin{equation}\label{lemconc}
\lim_{m\rightarrow\infty}\frac{I_{m-2}\left(\arccos\e^{-m^{\gamma}}
\right)}{I_{m-2}(\pi)}=\frac{1}{2}.
\end{equation}
An elementary proof of this fact can be found in Lemma \ref{apdx2} 
of Appendix~A. 

\begin{corollary}\label{polycor2}
For all $(m,n)$ such that $m\leq n$ let $A^{[m,n]}$ be 
a uniform random $n\times m$ matrix. Then 
\begin{equation*}
\limsup_{m\rightarrow\infty}\left(\sup_{n\geq m}
\frac{\bfE\left[\log\CN(A^{[m,n]})\right]}
{m}\right)\leq\log 2. 
\end{equation*}
\end{corollary}

\proof
Let $\epsilon>0$ be a small real number and $X$ a binomially 
distributed random variable $X\sim\mathbf{Bin}(m+k,(1-\epsilon)/2)$. 
If $\Phi$ denotes the cumulative distribution function of the 
standard normal distribution, then the central limit theorem shows that 
\begin{align}
&\binom{m+k}{m}\left(\frac{1+\epsilon}{2}\right)^{k}
=\left(\frac{2}{1-\epsilon}\right)^m\bfP\left[X=m\right]\nonumber\\
&=\left(\frac{2}{1-\epsilon}\right)^m
\bfP\left[\frac{X-\frac{(m+k)(1-\epsilon)}{2}}
{\frac{\sqrt{(m+k)(1-\epsilon^2)}}{2}}\in
\left[\frac{m-k-1+\epsilon(m+k)}{\sqrt{(m+k)(1-\epsilon^2)}},
\frac{m-k+1+\epsilon(m+k)}{\sqrt{(m+k)(1-\epsilon^2)}}\right)\right]
\nonumber\\
&\approx\left(\frac{2}{1-\epsilon}\right)^m\left(
\Phi\left(\frac{m-k+1+\epsilon(m+k)}{\sqrt{(m+k)(1-\epsilon^2)}}\right)-
\Phi\left(\frac{m-k-1+\epsilon(m+k)}{\sqrt{(m+k)(1-\epsilon^2)}}
\right)\right)\label{approximate equality}\\
&<\left(\frac{2}{1-\epsilon}\right)^m\cdot
\frac{2}{\sqrt{(m+k)(1-\epsilon^2)}}\cdot\nonumber\\
&\qquad\cdot
\frac{1}{\sqrt{2\pi}}\exp\left\{-\frac{1}{2}\min\left(
\frac{(m-k+1+\epsilon(m+k))^2}{(m+k)(1-\epsilon^2)},
\frac{(m-k-1+\epsilon(m+k))^2}{(m+k)(1-\epsilon^2)}\right)\right\}.
\nonumber
\end{align}
The approximate equality \eqref{approximate equality} becomes 
asymptotically exact. Therefore, there exists a number 
$m_{\epsilon}^1\in\N$ such that for all $m\geq m^1_{\epsilon}$ we have 
\begin{equation}\label{expo}
\begin{split}
&\binom{m+k}{m}\left(\frac{1+\epsilon}{2}\right)^{k}\leq
\exp\left\{m\log\frac{2}{1-\epsilon}+\frac{1}{2}\log\frac{2}{\pi}
\right.\\
&\,\left.
-\frac{1}{2}\log(m+k)-\frac{1}{2}\min\left(
\frac{(m-k+1+\epsilon(m+k))^2}{(m+k)(1-\epsilon^2)},
\frac{(m-k-1+\epsilon(m+k))^2}{(m+k)(1-\epsilon^2)}\right)\right\}.
\end{split}
\end{equation}
Moreover, \eqref{lemconc} shows that there exists a number 
$m^2_{\epsilon}\in\N$ such that for all $m\geq m^2_{\epsilon}$ we have
\begin{equation}\label{conc2}
\frac{I_{m-2}\left(\arccos\e^{-\sqrt{m}}
\right)}{I_{m-2}(\pi)}\leq\frac{1+\epsilon}{2}.
\end{equation}
Equations \eqref{round A}, \eqref{expo} and \eqref{conc2} show that for 
all $m\geq\max(m^1_{\epsilon},m^2_{\epsilon})$, $n=m+k\geq m$ and 
$t\geq\sqrt{m}$,  
\begin{align*}
&\bfP\left[\log\CN(A^{[m,n]})\geq t\right]\leq 
\binom{m+k}{m}2m^{\frac{5}{2}}\left(\frac{1+\epsilon}{2}\right)^{k}\e^{-t}\\
&\;\leq\exp\left\{m\log\frac{2}{1-\epsilon}+\frac{3}{2}\log 2 +
\frac{5}{2}\log m-t\right\}.
\end{align*}
Lemma~\ref{polylemma} implies that for the same parameters, 
\begin{equation*}
\bfE\left[\log\CN(A^{[m,n]})\right]\leq
\max\left(m\log\frac{2}{1-\epsilon}+\frac{3}{2}\log 2 +
\frac{5}{2}\log m, \sqrt{m}\right)+1.
\end{equation*}
Since $\epsilon$ was arbitrary, the claim follows.  
\eproof
\hspace{1cm}

The asymptotic linearity in $m$ of the bound on $\bfE[\log\CN(A)]$ 
derived in the above corollary is due to the appearance of a 
binomial term in \eqref{round A}. This is largely an artifact of 
our specific analysis, and if the hunch of Remark~\ref{suspicion} 
is true, then the asymptotic behaviour for $n\geq m\gg 1$ is given 
by 
\begin{equation}\label{gamma}
\bfE[\log\CN(A)]=\Oh(m^\gamma)
\end{equation}
for any arbitrarily small real exponent $\gamma>0$. However, it will not be 
logarithmically small in $m$, because of the concentration of measure 
phenomenon. Thus, if the hunch is true, then \eqref{gamma} describes
the asymptotic behaviour for arbitrary $n\geq m$ and $m\gg 1$. On 
the other hand, when $n\geq 5m$, we can actually prove that
\eqref{gamma} holds true.

\begin{corollary}\label{polycor3}
Let $\{A^{[m,n]}\}$ be the set of random matrices defined in Corollary 
\ref{polycor2}. Then for any real exponent $\gamma>0$, 
\begin{equation*}
\limsup_{m\rightarrow\infty}\left(\sup_{n\geq 5m}
\frac{\bfE\left[\log\CN(A^{[m,n]})\right]}
{m^{\gamma}}\right)=0,
\end{equation*}
that is, $\bfE[\log\CN(A^{[m,n]})]$ grows more slowly than any algebraic
function of $m$ when $n\geq 5m$. 
\end{corollary}

\proof
We need to consider \eqref{expo} again. For $k\geq 4m$ and $\epsilon$ 
small enough we have 
\begin{equation*}
m\log\frac{2}{1-\epsilon}<\frac{1}{2}
\min\left(
\frac{(m-k+1+\epsilon(m+k))^2}{(m+k)(1-\epsilon^2)},
\frac{(m-k-1+\epsilon(m+k))^2}{(m+k)(1-\epsilon^2)}\right),
\end{equation*}
and then 
\begin{equation}
\binom{m+k}{m}\left(\frac{1+\epsilon}{2}\right)^k\leq
\exp\left\{\frac{1}{2}\log\frac{2}{\pi}\right\}\label{conc55}
\end{equation}
for all $m\geq m^1_{\epsilon}$. Moreover, \eqref{lemconc} shows 
that there exists a number $m^3_{\epsilon}$ such that for all 
$m\geq m^3_{\epsilon}$ 
\begin{equation}\label{conc3}
\frac{I_{m-2}\left(\arccos\e^{-m^{\frac{\gamma}{2}}}
\right)}{I_{m-2}(\pi)}\leq\frac{1+\epsilon}{2}.
\end{equation}
Equations \eqref{round A}, \eqref{conc2}, \eqref{conc55} and 
\eqref{conc3} together imply that 
\begin{equation*}
\bfP\left[\log\CN(A^{[m,n]})\geq t\right]\leq
\exp\left\{\frac{3}{2}\log2-\frac{1}{2}\log\pi+\frac{5}{2}\log m -t\right\}
\end{equation*}
for all $n\geq 5m$, $m\geq m^{\frac{\gamma}{2}}_{\epsilon}=
\max(m^1_{\epsilon},m^3_{\epsilon})$ and $t\geq
m^{\frac{\gamma}{2}}$. Finally, applying Lemma~\ref{polylemma}, we get 
\begin{equation*}
\bfE\left[\log\CN(A^{[m,n]})\right]\leq 1+\max\left(
\frac{3}{2}\log2-\frac{1}{2}\log\pi+\frac{5}{2}\log m,m^{\frac{\gamma}
{2}}\right).
\end{equation*}
Dividing by $m^{\gamma}$ and taking limits, the result follows. 
\eproof

\section{A Final Remark About the Case $n<m$}

The development in the previous sections assumes $n\geq m$. The case 
$n<m$ has been dealt with in~\cite{ChC01}, where it is proved that 
$$
  \bfE[\log\CN(A)]\leq \frac{5}{2}\log n+2.
$$
\medskip

\bibliographystyle{plain}

\section*{Appendix A: A Concentration of Measure
Inequality} 

Let $X\sim\UN(\sphere^{m-1})$ and $p\in\sphere^{m-1}$. 
Note that for $\rho<\pi/2$, 
\begin{align}
\bfP&\left[X\in\Cap(p,\rho)\right]=
\frac{I_{m-2}(\rho)}{I_{m-2}(\pi)}\nonumber\\
&=\frac{1}{2}
\frac{\int_{0}^{\rho}\sin^{m-2}\tau d\tau}{\int_0^{\rho}
\sin^{m-2}\tau d\tau +\int_{\rho}^{\frac{\pi}{2}}\sin^{m-2}
\tau d\tau}\label{round 1}\\
&=\frac{1}{2}\frac{\int_0^{\rho}\sin^{m-2}\tau d\tau}
{\int_{\rho}^{\frac{\pi}{2}}\sin^{m-2}\tau d\tau}
\left(1+\frac{\int_0^{\rho}\sin^{m-2}\tau d\tau}{\int_{\rho}
^{\frac{\pi}{2}}\sin^{m-2}\tau d\tau}\right)^{-1}\nonumber\\
&=\frac{1}{2}\left(\xi-\xi^2+\xi^3-\dots\right)
=\frac{1}{2}\xi+\Oh(\xi^2),\nonumber
\end{align}
where the last line holds if 
\begin{equation*}
\xi:=\frac{\int_{0}^{\rho}\sin^{m-2}\tau d\tau}{\int_{\rho}
^{\frac{\pi}{2}}\sin^{m-2}\tau d\tau}<1.
\end{equation*}
Now, for $\rho\ll \pi/2$ we have $\xi=\Oh(\rho^{m-1})$, and 
hence, $\bfP[X\in\Cap(p,\rho)]=\Oh(\rho^{m-1})$ as one would 
expect. Likewise, one expects intuitively that if 
\begin{equation}\label{sharper meaning}
\frac{\frac{\pi}{2}-\rho}{\frac{\pi}{2}}\ll 1
\end{equation}
then 
\begin{equation}\label{round 2}
\bfP\left[X\in\Cap(p,\rho)\right]=\frac{1}{2}-\Oh\left(
\frac{\frac{\pi}{2}-\rho}{\frac{\pi}{2}}\right),
\end{equation}
and this is indeed the case. However, it is somewhat surprising 
that for large $m$, the expression 
\begin{equation}\label{size}
\frac{\frac{\pi}{2}-\rho}{\frac{\pi}{2}} 
\end{equation}
has to be extremely small indeed before the order \eqref{round 2} is 
observed. In fact, if \eqref{size} decreases to zero at an algebraic 
rate as a function of the dimension $m$, then $\bfP\left[
X\in\Cap(p,\rho)\right]$ converges to 
zero. We are not going to prove this property here, although an 
elementary proof can be given along the lines of Lemma \ref{apdx2} 
below, but we remark that this is a special case of a type of 
properties of high-dimensional probability distributions that are 
jointly referred to as the {\em concentration of measure phenomenon}. 
See e.g.\ \cite{Ledoux} for a good account of this theory. 
The purpose of this appendix is in some sense to get around the 
adverse effects of the concentration of measure phenomenon and to 
show that if the expression \eqref{size} is exponentially small in
terms of $m$ then \eqref{round 2} is asymptotically observed. In fact, 
we are going to prove a slightly weaker result which is sufficient for 
the purposes of the analysis of Section \ref{mgg1}.

\begin{lemma}\label{apdx2}
Let $\gamma>0$ be a constant, $p\in\sphere^{m-1}$ and 
$X\sim\UN(\sphere^{m-1})$. Then 
\begin{equation*}
\lim_{m\rightarrow\infty}\bfP\left[
X\in\Cap\left(p,\arccos\e^{-m^{\gamma}}\right)\right]=\frac{1}{2}.
\end{equation*}
\end{lemma}

\proof
Let $\theta\in(0,1/2)$ and let $m_{\theta}\in\N$ be such that 
\begin{equation*}
\theta<\frac{1-\frac{1}{m_{\theta}-3}}{2}. 
\end{equation*}
Then, for $m\geq m_{\theta}$ equation \eqref{round 1} implies that 
\begin{equation}\label{round 4}
\bfP[X\in\Cap(p,\rho)]<\theta\Leftrightarrow
\int_{0}^{\rho}\sin^{m-2}\tau d\tau<2\theta\int_0^{\frac{\pi}{2}}
\sin^{m-2}\tau d\tau.
\end{equation}
It is easy to show by induction and partial integration that 
\begin{equation}\label{round 5}
\int_{\rho}^{\frac{\pi}{2}}\sin^{m}\tau d\tau=\begin{cases}
&\frac{1}{m}(\cos\rho)\left[\sin^{m-1}\rho+\sum_{k=0}^{\frac
{m-3}{2}}\frac{(m-1)(m-3)\dots (m-2k-1)}{(m-2)(m-4)\dots 
(m-2k-2)}\sin^{m-2k-3}\rho\right]\\
&\hspace{2cm}\text{ if } m\text{ is odd},\\
&\frac{1}{m}(\cos\rho)\left[\sin^{m-1}\rho+\sum_{k=0}^{\frac
{m}{2}-2}\frac{(m-1)(m-3)\dots (m-2k-1)}{(m-2)(m-4)\dots 
(m-2k-2)}\sin^{m-2k-3}\rho\right]\\
&\quad+\frac{(m-1)(m-3)\dots 3\cdot 1}
{m(m-2)\dots 2}(\frac{\pi}{2}-\rho)\text{ if } m\text{ is even}.
\end{cases}
\end{equation}
In particular,
\begin{equation}\label{round 6}  
\int_{0}^{\frac{\pi}{2}}\sin^m\tau d\tau=\begin{cases}
\frac{(m-1)(m-3)\dots 4\cdot 2}{m(m-2)\dots 3\cdot 1}
&\text{ if }m \text{ is odd},\\
\frac{(m-1)(m-3)\dots 3\cdot 1}{m(m-2)\dots 2}\cdot\frac{\pi}{2}
&\text{ if }m\text{ is even}.
\end{cases}
\end{equation}
It follows from \eqref{round 5} and \eqref{round 6} that 
\begin{align}
\int_{\rho}^{\frac{\pi}{2}}\sin^m\tau d\tau&<\begin{cases}
&\frac{(m-1)(m-3)\dots 2}{m(m-2)\dots 1}\cdot\frac{1-sin^{m+1}\rho}
{\cos\rho}=\frac{1-\sin^{m+1}\rho}{\cos\rho}\int_0^{\frac{\pi}{2}}
\sin^m\tau d\tau\\
&\hspace{2cm}\text{if }m\text{ is odd},\\
&\\
&\frac{(m-1)(m-3)\dots 3}{m(m-2)\dots 2}\left[
(\tan\rho)(1-\sin^m\rho)+\frac{\pi}{2}-\rho\right]\\
&\quad=\left[\frac{2\sin\rho(1-\sin^m\rho)}{\pi\cos\rho}+1-\frac{2\rho}{\pi}
\right]\cdot\int_{0}^{\frac{\pi}{2}}\sin^m\tau d\tau\\
&\hspace{2cm}\text{ if }m\text{ is even}.
\end{cases}\nonumber\\
&<\frac{(1-\sin^{m+1}\rho)+2(1-\sin^2\rho)}{\cos\rho}\cdot
\left(\int_0^{\frac{\pi}{2}}\sin^m\tau d\tau\right)\qquad
\text{(in both cases)},\nonumber
\end{align}
where the last inequality holds at least for $\pi/4\leq\rho<\pi/2$.
Therefore, 
\begin{align}\label{round 8}
\int_{0}^{\rho}\sin^m\tau d\tau&=\int_{0}^{\frac{\pi}{2}}\sin^m\tau d\tau
-\int_{\rho}^{\frac{\pi}{2}}\sin^m\tau d\tau\nonumber\\
&>\left[1-\frac{(1-\sin^{m+1}\rho)+2(1-\sin^2\rho)}{\cos\rho}\right]\cdot
\left(\int_0^{\frac{\pi}{2}}\sin^m\tau d\tau\right).
\end{align}
Note that the combination of \eqref{round 4} and \eqref{round 8} 
implies that for $\pi/4\leq\rho<\pi/2$, 
\begin{equation}\label{see?}
2\theta<1-\frac{(1-\sin^{m+1}\rho)+2(1-\sin^2\rho)}{\cos\rho}
\Rightarrow\bfP\left[X\in\Cap(p,\rho)\right]\geq\theta.
\end{equation}

Now let the sequence $(\rho_m)_{\N}$ be defined by 
\begin{equation*}
\rho_m=\arccos\e^{-m^{\gamma}},
\end{equation*}
and note that 
\begin{equation*}
1-\frac{(1-\sin^{m+1}\rho_m)+2(1-\sin^2\rho_m)}{\cos\rho_m}
=1-\e^{m^{\gamma}(1-2(m+1))}-2\e^{-m^{\gamma}}\stackrel
{m\rightarrow\infty}{\longrightarrow}1.
\end{equation*}
Therefore, for $m\geq m_{\theta}$ large enough,
$\rho_m\in(\pi/4,\pi/2)$ and the condition on the left hand side of 
\eqref{see?} is satisfied. This shows that 
\begin{equation*}
\lim_{m\rightarrow\infty}\bfP\left[
X\in\Cap\left(p,\arccos\e^{-m^{\gamma}}\right)\right]\geq\theta,
\end{equation*}
and since this is true for any $\theta\in(0,1/2)$, this proves that 
$\lim_{m\rightarrow\infty}\bfP\left[
X\in\Cap\left(p,\arccos\e^{-m^{\gamma}}\right)\right]\geq\frac{1}{2}$.
Moreover, the inequality $\lim_{m\rightarrow\infty}\bfP\left[
X\in\Cap\left(p,\arccos\e^{-m^{\gamma}}\right)\right]\leq\frac{1}{2}$
is trivial, and the result follows.
\eproof

\section{Appendix B: Complexity of the Relaxation Method}

The purpose of this appendix is to make a convincing argument that 
the complexity of relaxation methods for the solution of the linear 
system $Ax\leq 0,\;x\neq 0$ is proportional to $\CN(A)^2$. In a sense, 
this fact is in the general knowledge of researchers familiar with both 
the relaxation method and condition numbers, as conversations 
with Marina Epelman, Rob Freund and Dan Spielman confirmed. Moreover, 
this fact has been implicitly stated in the relaxation method
literature for decades, and all it takes to make it explicit is to 
translate well-known results from the relaxation method literature 
into the language of condition numbers. For lack of an explicit 
reference, let us give such an example here. 

The algorithm we consider is the so-called perceptron algorithm 
\cite{Rosenblatt}. When applied to solving a strictly feasible 
system $Ax<0$, where $A\in\R^{n\times m}$ with unit row vectors, 
this algorithm starts from an initial point 
$x_0\in\R^m$ and constructs an iterative sequence of points 
$(x_i)_{\N}\subset\R^m$ as follows: if $Ax_i<0$ then $x_i$ is a 
solution and the algorithm stops. Otherwise, a row vector $a^{[i]}$ 
from $A$ is chosen so that $a^{[i]}x_i\geq 0$, and the next point 
is computed by $x_{i+1}=x_i-a_i$. 

The usual convergence analysis then proceeds as follows, see e.g.\ 
\cite{Block}: let $w\in\R^m$ be a solution of $Ax<0$, let 
$\alpha=\max_{1\leq j\leq n}a_j w$, where $a_j$ is the $j$-th 
row vector of $A$, and let $w^*=w/|\alpha|$. It is then easy to 
show that if the algorithm has not stopped before or during iteration 
$i$, that is if $x_{i+1}\neq x_i$, then  
\begin{equation*}
\|x_{i+1}-w^*\|^2\leq\|x_i-w^*\|^2-1.
\end{equation*}
This shows that at most $\|x_0-w^*\|^2$ iterations can take place. 
For simplicity, we can choose $x_0$ to be the origin, so that the 
algorithm has complexity $\Oh(\|w^*\|^2)$. 

Now note that $\alpha=\cos\theta(A,w)$, where we use the notation of 
Section \ref{notions}. Without loss of generality we may assume that 
$w$ is a unit vector, so that $\|w^*\|=|\cos\theta(A,w)|^{-1}$. In 
order to minimise the complexity estimate, we need to choose
$w=\bar{x}$, so that we find that the algorithm terminates after at 
most $\CN(A)^2$ iterations. 

\end{document}